\documentclass{amsart}
\usepackage{amsfonts}
\usepackage{amssymb}
\usepackage{amsthm}
\usepackage{eucal}
\usepackage[all]{xy}
\usepackage{graphicx}
\usepackage[ansinew]{inputenc}
\newtheorem{thm}{Theorem}[section]
\newtheorem{cor}[thm]{Corollary}
\newtheorem{lem}[thm]{Lemma}
\newtheorem{prop}[thm]{Proposition}
\theoremstyle{definition}

\theoremstyle{remark}
\newtheorem{rem}[thm]{Remark}
\numberwithin{equation}{section}

\newcommand{\set}[1]{\left\{#1\right\}}

\newcommand{\To}{\longrightarrow}

\newcommand{\mo}{\mathrm{\mathbf{mod}}}
\newcommand{\sub}{\mathrm{\mathbf{sub}}}
\newcommand{\rep}{\mathrm{\mathbf{rep}}}
\newcommand{\fp}{\mathrm{\mathbf{fp}}}

\newcommand{\Ab}{\mathrm{\mathbf{Ab}}}
\newcommand{\M}{\mathrm{\mathbf{M}}}
\newcommand{\pair}{\mathrm{\mathbf{pair}}}
\newcommand{\n}[1]{\textsl{\textbf{#1}}}

\def\ker{\operatorname{Ker}}

\def\coker{\operatorname{Coker}}
\def\card{\operatorname{card}}
\def\hom{\operatorname{Hom}}
\def\ext{\operatorname{Ext}}
\def\endo{\operatorname{End}}
\def\colim{\mathop{\operatorname{colim}}}
\newcommand{\grupo}[1]{\langle #1\rangle}

\def\Z{{\mathchoice {\hbox{$\mathsf{\textstyle{Z\kern-0.4em Z}}$}}
{\hbox{$\mathsf{\textstyle{ Z\kern-0.4em Z}}$}}
{\hbox{$\mathsf{\scriptstyle{ Z\kern-0.3em Z}}$}}
{\hbox{$\mathsf{\scriptscriptstyle{ Z\kern-0.2em Z}}$}}}}
\def\N{{\mathrm{I\!N}}}
\def\st{\stackrel} 
\def\ul{\underline} 
\def\ol{\bar} 
\newcommand{\um}[1]{\mathcal{#1}} 
\newcommand{\uf}[1]{\mathbb{#1}} 
\newcommand{\mat}[1]{\mathsf{#1}} 

\numberwithin{equation}{section}
\renewcommand{\theequation}{\thesection.\Alph{equation}}

\def\NN{\N_0}

\begin{document}


\title[Representation theory of some infinite-dimensional algebras\dots]
{Representation theory of some infinite-dimensional algebras arising in continuously controlled algebra and topology}%
\author{Fernando Muro}%
\address{Departamento de Geometr\'{\i}a y Topolog\'{\i}a. Universidad de Sevilla.}%
\email{fmuro@us.es}%



\maketitle

\tableofcontents
\section{Introduction}\label{intro}

Suppose that the discrete set $\NN$ of non-negative integers is
embedded $\NN\subset X$ in a compact metrizable space $X$, and let
$E=\NN'\subset X$ be the derived set (i. e. points of $X$ which
contain infinitely many points of $\NN$ in any neighborhood).
Consider the set $R(E)$ of $\NN\times\NN$ matrices
$(\mat{a}_{ij})_{i,j\in\NN}$ with entries in a (unital
associative) ring $R$ such that if $\set{i_n}_{n\in\NN},
\set{j_n}_{n\in\NN}\subset\NN$ are sequences convergent in $X$ to
different points then the vector $(\mat{a}_{i_nj_n})_{n\in\NN}$ is
almost all zero. This set is an $R$-algebra with the usual matrix
operations. Any compact metrizable space can arise as $E$ in this
way. In fact the isomorphism class of the algebra $R(E)$ only
depends on $E$. These algebras are Morita equivalent to some
additive categories of free $R$-modules continuously controlled at
infinity by $E$ appearing in the literature. These categories play
an important role in many areas such as controlled homotopy
theory, proper homotopy theory, $C^*$-algebra theory, $K$-theory
and $L$-theory, see for example \cite{quinn}, \cite{canc},
\cite{ccstar}, \cite{iht} and \cite{projdim}.

The elementary properties of the algebras $R(E)$ have been studied
by Baues-Quintero (\cite{iht}) for $R=\Z$ the integers. If $E$ is
zero-dimensional, $R(E)$ is a particular case of the rings
considered by Farrell-Wagoner (\cite{fw}). When $E=*$ is a
singleton $R(E)=\mathrm{RCFM}(R)$ is the well-known algebra of
row-column-finite (or locally finite) matrices over $R$. This
algebra has been studied from a purely ring-theoretical point of
view, see for example \cite{rcfmurcia}. It was also used by
Wagoner (\cite{deloopw}) to construct deloopings in algebraic
$K$-theory.

In this paper we concentrate on the representation theory of the
algebras $k(E)$, where $k$ is any field.

Representation theory considers the \emph{decomposition problem}
in a small additive category $\mathbf{A}$. A solution to this
problem consists of a set of objects (which we call
\emph{elementary objects}) and of a set of isomorphisms
(\emph{elementary isomorphisms}) between finite direct sums of
elementary objects. These sets must satisfy the following
properties: any object in $\mathbf{A}$ is isomorphic to a finite
direct sum of elementary ones, and any isomorphism relation
between two such direct sums
can be derived from the elementary isomorphisms. 
Notice that this is exactly a presentation of the abelian monoid
$\operatorname{Iso}(\mathbf{A})$ of isomorphisms classes of
objects in $\mathbf{A}$. The trivial solution is taking all
objects as elementary objects and all isomorphisms as elementary
isomorphisms, however one is often interested in solutions
minimizing the cardinal of the sets of elementary objects and
isomorphisms.

We say that $\mathbf{A}$ has \emph{finite representation type} if
there exists a finite set of elementary objects, or equivalently
$\operatorname{Iso}(\mathbf{A})$ is finitely presented. The
representation type of $\mathbf{A}$ is \emph{wild} if a solution
to the decomposition problem in $\mathbf{A}$ would yield a
solution to the decomposition problem in the category of
finite-dimensional modules over a polynomial $k$-algebra in two
non-commuting variables. Otherwise $\mathbf{A}$ has \emph{tame}
representation type. If $\mathbf{A}$ has wild representation type
the word problem for finitely presented groups, which is
undecidable, can be embedded in the decomposition problem in
$\mathbf{A}$, hence one can not expect to get satisfactory
solutions in this case. The representation type of an algebra $A$
is that of the category $\fp(A)$ of finitely presented (right)
$A$-modules.

One of the main results of this paper is the following theorem,
where we compute the representation type of the algebra $k(E)$ in
terms of the cardinal of $E$, without restrictions on the ground
field $k$.

\begin{thm}\label{rep}
The representation type of $k(E)$ is
\medskip
\begin{center}
\begin{tabular}{|c|c|}
  \hline
  $\card E$ & type \\
  \hline
  $<4$ & finite \\
  $=4$ & tame \\
  $>4$ & wild \\
  \hline
\end{tabular}
\end{center}
\end{thm}

\bigskip


In the finite and tame cases we construct explicit presentations
of $\operatorname{Iso}(\fp(k(E)))$. Moreover, for $E$ finite, we
prove that there are presentations of
$\operatorname{Iso}(\fp(k(E)))$ with a finite number of elementary
isomorphisms and we compute them. These presentations satisfy the
next properties.

\begin{thm}\label{generel}
If $\card E$ is finite there are solutions to the decomposition
problem in $\fp(k(E))$ with the next cardinals of elementary
modules and isomorphisms
\medskip
\begin{center}
\begin{tabular}{|c|c|c|}
  \hline
  $\card E$ & modules & isomorphisms \\
  \hline
  $1$ & $6$ & $6$ \\
  $2$ & $12$ & $12$ \\
  $3$ & $21$ & $18$ \\
  $\geq 4$ & $\geq\aleph_0$ & $6\card E$ \\
  \hline
\end{tabular}
\end{center}
\end{thm}





There are two key steps in the proof of these results. The first
is to solve the decomposition problem for finitely presented
$\mathrm{RCFM}(k)$-modules. The second is to relate the
decomposition problem in $\fp(k(E))$ when $\card E=n$ is finite to
the decomposition problems in
$\fp(\mathrm{RCFM}(k))$ and in the category of finite-dimensional $n$-subspaces.

We shall use the pro-category of pro-vector spaces and the inverse
limit functor to construct invariants detecting isomorphism types
of finitely presented $k(E)$-modules. Moreover, we shall not
usually work directly with the categories of finitely presented
$k(E)$-modules but with the equivalent categories of finitely
presented modules over certain small additive categories. This
setting allows more flexibility and technical proofs become less
complicated than if we use $k(E)$-modules.

In order to facilitate the reading we now describe the contents of
this paper. In the next section we briefly recall from \cite{rso}
the basic tools of ringoid theory that we need. Afterwards, in
Section \ref{unosring}, we introduce the ringoids which are Morita
equivalent to the algebras $R(E)$ and establish their basic
properties. For this we use the approach in \cite{iht},
generalizing some results in this reference for $R=\Z$ to
arbitrary rings. We put emphasis on the case $E$ finite because we
shall always work under this assumption (even for the proof of
Theorem \ref{rep}, see Remark \ref{reduc}). In Section \ref{count}
we construct an embedding of the decomposition problem for
countably presented $R$-modules into the decomposition problem for
finitely presented $R(E)$-modules, where $E$ is any non-empty
compact metrizable space. Section \ref{repro} contains basic facts
about pro-categories. In Section \ref{nume} we construct some
invariants of the isomorphism class of a finitely presented
$k(E)$-module, $E$ finite. These invariants are used in Section
\ref{rcfm} to classify finitely presented $k(E)$-modules when
$E=*$ is just one point and $k(E)=\mathrm{RCFM}(k)$. In particular
we prove that this $k$-algebra has finite representation type. The
classification theorem (Theorem \ref{clasifR}) is derived from
several technical lemmas. In Section \ref{n} we recall the
definition and representation theory of the $n$-subspace quiver.
We also define the class of rigid $n$-subspaces, which plays an
important role in what follows. We show that all but $3n$
indecomposable representations of the $n$-subspace quiver are
rigid $n$-subspaces. In Section \ref{modandn} we relate the
representation theories of both $k(E)$ and the $n$-subspace
quiver, where $n$ is the cardinal of $E$. The properties of this
relation are established through many technical results which lead
us to complete the proofs of Theorems \ref{rep} and \ref{generel}
in Section \ref{oju}. In this last section we compute the
structure of the monoid $\operatorname{Iso}(\fp(k(E)))$ (Theorem
\ref{jander}) and give a classification theorem for finitely
presented $k(E)$-modules (Corollary \ref{clasclas}) for $E$
finite. This classification theorem explicitly describes the
(finite) set of elementary isomorphisms, and also the set of
elementary objects when $E$ has less than $5$ points. We include
an Appendix with some computations of $\ext^1$ groups of finitely
presented $k(E)$-modules, $E$ finite. These computations will be
applied to proper homotopy theory in a forthcoming paper
(\cite{propera2n}).

\subsection{Notation and conventions} In this paper all rings and
algebras are associative with unit. We use bold letters
$\mathbf{C}$ for categories, $R$ for an arbitrary
(non-commutative) ring, $\Z$ for the ring of integers, and $k$ for
fields. As usual $\N=\set{1,2,3,\dots}$ is the set of natural
numbers, and $\NN=\N\cup\set{0}$ the free abelian monoid with one
generator.

Capital sans serif letters $\mat{A}$ are names of matrices with
entries in some ring. Here all matrices are square matrices
indexed by $\NN$, and the entry of a matrix $\mat{A}$
corresponding to the subindexes $i,j\in\NN$ is denoted by
$\mat{a}_{ij}$, it is $\mat{A}=(\mat{a}_{ij})_{i,j\in\NN}$. The
identity matrix is denoted by $\mat{I}$, it is defined as
$\mat{i}_{ii}=1$ $(i\in\NN)$ and $\mat{i}_{ij}=0$ for $i\neq j$.
The entries of the transposed matrix
$\mat{A}^t=(\mat{a}^t_{ij})_{i,j\in\NN}$ of $\mat{A}$ are
$\mat{a}_{ij}^t=\mat{a}_{ji}$ $(i,j\in\NN)$. Vectors are denoted
by $(\mat{v}_i)_{i=1}^n$ or $(\mat{v}_n)_{n\in\NN}$ provided they
have a finite or an infinite countable number of entries. We
regard vectors as column matrices, hence matrices act on vectors
on the left.

\section{Ringoids and modules}\label{rym}

A \emph{ringoid} $\mathbf{R}$ is a category whose morphism sets
$\hom_\mathbf{R}(X,Y)$ are abelian groups in such a way that
composition is bilinear. The endomorphism set
$\endo_\mathbf{R}(X)=\hom_\mathbf{R}(X,X)$ of an object $X$ has a
ring structure with product given by composition of morphisms.
Conversely any ring $R$ is identified with the ringoid with a
single object whose endomorphism set is $R$. An \emph{additive
category} is a ringoid with finite biproducts (direct sums). In
this section we recall basic facts about modules over a small
ringoid. Our main reference for this subject is \cite{rso}.

An \emph{additive functor} between ringoids is a functor which
induces homomorphisms between morphism sets. Let $\Ab$ be the
category of abelian groups. A \emph{right-$\mathbf{R}$-module}
$\um{M}$ is an additive functor
$\um{M}\colon\mathbf{R}^{op}\rightarrow\Ab$. Morphisms of
right-$\mathbf{R}$-modules are natural transformations, and the
category of right-$\mathbf{R}$-modules is denoted by
$\mo(\mathbf{R})$ whenever $\mathbf{R}$ is small.
Left-$\mathbf{R}$-modules are the same thing as
right-$\mathbf{R}^{op}$-modules, where $\mathbf{R}^{op}$ is the
opposite category, so every statement about right-modules has a
convenient translation to left-modules. From now on every module
is a right-module unless we state the contrary.

There is a Yoneda full inclusion of categories $\mathbf{R}\subset
\mo(\mathbf{R})$ which sends an object $X$ in $\mathbf{R}$ to the
associated contravariant representable functor
$\hom_\mathbf{R}(-,X)$. These $\mathbf{R}$-modules are said to be
\emph{finitely generated free}. They are projective by Yoneda's
lemma.

An $\mathbf{R}$-module $\um{M}$ is \emph{finitely presented} (f.
p.) if it is the cokernel of a morphism between two finite direct
sums of finitely generated free $\mathbf{R}$-modules. The cokernel
of a morphism between f. p. modules is also f. p. In particular
direct summands of f. p. modules are f. p. We write
$\fp(\mathbf{R})\subset\mo(\mathbf{R})$ for the full subcategory
of f. p. $\mathbf{R}$-modules. If $\mathbf{R}=\mathbf{A}$ is an
additive category, then a f. p. $\mathbf{A}$-module $\um{M}$ is in
fact the cokernel in $\mo(\mathbf{A})$ of a morphism
$\varphi\colon {X}_1\rightarrow {X}_0$ in $\mathbf{A}$. One can
readily check that if $\um{N}=\coker[\psi\colon {Y}_1\rightarrow
{Y}_0]$ is another f. p. $\mathbf{A}$-module, any morphism
$\tau\colon \um{M}\rightarrow \um{N}$ is represented by a morphism
$\tau_0\colon {X}_0\rightarrow {Y}_0$ such that there exists
$\tau_1\colon {X}_1\rightarrow {Y}_1$ with
$\tau_0\varphi=\psi\tau_1$. Another morphism $\tau_0'\colon \colon
{X}_0\rightarrow {Y}_0$ represents $\tau$ if and only if there
exists $\eta\colon {X}_0\rightarrow {Y}_1$ with
$\tau_0+\psi\eta=\tau_0'$. More precisely, let $\pair(\mathbf{A})$
be the additive category whose objects are morphisms in
$\mathbf{A}$, and morphisms
$\tau=(\tau_1,\tau_0)\colon\varphi\rightarrow\psi$ are commutative
squares
$$\xymatrix{{X}_1\ar[r]^\varphi\ar[d]_{\tau_1}&{X}_0\ar[d]^{\tau_0}\\
{Y}_1\ar[r]^\psi&{Y}_0}$$ There is an obvious functor
$\coker\colon\pair(\mathbf{A})\rightarrow\fp(\mathbf{A})$ given by
taking cokernels. We define in $\pair(\mathbf{A})$ the natural
equivalence relation $\sim$ with $\tau\sim\tau'$ if there exists
$\eta\colon {X}_0\rightarrow {Y}_1$ satisfying
$\tau_0+\psi\eta=\tau_0'$.

\begin{prop}\label{cokereq}
The functor $\coker$ factors through the quotient category
$\pair(\mathbf{A})/\!\!\sim$ and the induced functor
$\coker\colon\pair(\mathbf{A})/\!\!\sim\;
\rightarrow\fp(\mathbf{A})$ is an equivalence of categories.
\end{prop}

Any additive functor $\uf{F}\colon\mathbf{R}\rightarrow\mathbf{S}$
between small ringoids induces two ``change of coefficients''
additive functors
$\uf{F}^*\colon\mo(\mathbf{S})\rightarrow\mo(\mathbf{R})$ and
$\uf{F}_*\colon\mo(\mathbf{R})\rightarrow\mo(\mathbf{S})$. The
first one is exact and sends an $\mathbf{S}$-module $\um{M}$ to
the composite $\uf{F}^*\um{M}=\um{M}\uf{F}$. The second one is
left-adjoint to $\uf{F}^*$ ($\uf{F}_*$ is the left additive Kan
extension along $\uf{F}$, see \cite{rso} 6) and hence right-exact.
Moreover, the next diagram commutes
\begin{equation}\label{conmuta}
\xymatrix{\mathbf{R}\ar[r]^{\uf{F}}\ar[d]_{\mathrm{Yoneda}}&\mathbf{S}\ar[d]^{\mathrm{Yoneda}}\\
\mo(\mathbf{R})\ar[r]^{\uf{F}_*}&\mo(\mathbf{S})} \end{equation}
In addition if $\uf{F}$ is full and faithful then so is
$\uf{F}_*$, and in this case $\uf{F}^*\uf{F}_*$ is naturally
equivalent to the identity, see \cite{borceux1} 3.4.1. This
follows from the fact that any $\mathbf{R}$-module admits a
projective resolution by (arbitrary) direct sums of finitely
generated free ones. The functor $\uf{F}_*$ restricts to the full
subcategories of f. p. modules.

If we identify the endomorphism ring ${\endo}_\mathbf{R}(X)$ of an
object $X$ with the full subcategory of $\mathbf{R}$ whose unique
object is $X$ the change of coefficients $\uf{F}^*$ induced by the
inclusion $\uf{F}\colon\endo_\mathbf{R}(X)\subset\mathbf{R}$ is
the evaluation functor
$$ev_X=\uf{F}^*\colon\mo(\mathbf{R})\To\mo({\endo}_\mathbf{R}(X))\colon \um{M}\mapsto \um{M}(X).$$
The next proposition is an useful criterion to detect when a
ringoid is Morita equivalent to a ring. It is a consequence of
\cite{rso} 8.1.

\begin{prop}\label{mitchell}
If every object in $\mathbf{R}$ is a retract of $X$ then the
evaluation functor $ev_X$ is an additive equivalence of abelian
categories which restricts to an equivalence between the full
subcategories of finitely presented modules.
\end{prop}

In a more categorical language ringoids are defined as categories
enriched over the monoidal category of abelian groups with the
usual tensor product, compare \cite{borceux2} 6.2. For any ring
$R$ one can consider the monoidal category of $R$-$R$-bimodules
with the $R$-tensor product and define an \emph{$R$-ringoid} as a
category enriched over it. This is the same as an $R$-category in
the sense of \cite{rso} when $R$ is commutative. If $\mathbf{R}$
is an $R$-ringoid the endomorphism ring $\endo_\mathbf{R}(X)$ of
an object $X$ is in fact an $R$-algebra. In this case
$\mathbf{R}$-modules and morphisms between them take values in the
category of (right) $R$-modules in a natural way.

\section{The algebras $R(E)$ and related additive
categories}\label{unosring}

Given a ring $R$ and a set $A$ we write $R\grupo{A}$ for the free
$R$-module with basis set $A$. Free $R$-modules are
$R$-$R$-bimodules, hence the additive category of free $R$-modules
and right-$R$-module homomorphisms is an $R$-ringoid. The
\emph{carrier} of an element $x\in R\grupo{A}$ is the (finite) set
$\operatorname{carr}(x)\subset A$ such that
$z\in\operatorname{carr}(x)$ if $z$ appears with a non-trivial
coefficient in the linear expansion of $x$.

For every non-empty compact metrizable space $E$ there exists
another one $X$ containing $E$ such that the complement $Y=X-E$ is
dense in $X$. The triple $\ol{T}=(X,Y,E)$ can always be chosen to
be a tree-like space in the sense of \cite{iht} III.1.1. But one
can also take $X$ to be the (unreduced) cone over $E$,
$X=CE=E\times[0,1]/E\times\set{1}$. Here we identify $E$ with
$E\times\set{0}$ inside the cone $CE$.

A \emph{free $\ol{T}$-controlled $R$-module} $R\grupo{A}_\alpha$
is a free $R$-module $R\grupo{A}$ together with a function
$\alpha\colon A\rightarrow Y$, called \emph{height function}, such
that $\alpha^{-1}(K)$ is finite for every compact subspace
$K\subset Y$. The set $A$ is necessarily countable and the derived
set of $\alpha(A)$ in $X$ satisfies $\alpha(A)'\subset E$. This
derived set is the \emph{support} of $R\grupo{A}_\alpha$.
\emph{Controlled homomorphisms} $\varphi\colon
R\grupo{A}_\alpha\rightarrow R\grupo{B}_\beta$ are homomorphisms
between the underlying $R$-modules such that for every $x\in E$
and every neighborhood $U$ of $x$ in $X$ there exists another
neighborhood $V\subset U$ of $x$ in $X$ such that if $a\in A$
satisfies $\alpha(a)\in V$ then
$\beta(\operatorname{carr}(\varphi(a)))\subset U$. The category
$\M_R(\ol{T})$ of free $\ol{T}$-controlled $R$-modules and
controlled homomorphisms is a small additive category. Moreover,
it is an $R$-ringoid. The sum and $R$-actions on morphism sets are
given by those of the underlying free $R$-module homomorphisms,
and the direct sum of two objects is $R\grupo{A}_\alpha\oplus
R\grupo{B}_\beta=R\grupo{A\sqcup B}_{(\alpha,\beta)}$, where
$A\sqcup B$ is the disjoin union of sets and $(\alpha,\beta)\colon
A\sqcup B\rightarrow Y$ is defined as $\alpha$ over $A$ and
$\beta$ over $B$. 


\begin{rem}\label{zzz}
The category $\M_R(\ol{T})$ is defined in \cite{iht} III.4.7 for
$\ol{T}$ a tree-like space. However, as it is pointed out in the
Remark after that definition, it is equivalent to the category
$\mathcal{B}(X,E;R)$ in \cite{canc}. In particular $\M_R(\ol{T})$
only depends on $E$ up to equivalence of categories preserving
supports of objects (in fact equivalence of $R$-ringoids), see
1.23 and 1.24 in \cite{canc}.
\end{rem}

The next proposition shows that free $\ol{T}$-controlled
$R$-modules are classified by the underlying $R$-module and the
support.

\begin{prop}\label{clasiiso}
Two free $\ol{T}$-controlled $R$-modules $R\grupo{A}_\alpha$,
$R\grupo{B}_\beta$ are isomorphic if and only if the next two
conditions are satisfied:
\begin{enumerate}
\item The underlying $R$-modules are isomorphic $R\grupo{A}\simeq
R\grupo{B}$,

\item both have the same support $\alpha(A)'=\beta(B)'$.
\end{enumerate}
If the supports are non-empty then condition $(1)$ is
automatically satisfied. Furthermore, any compact subset $K\subset
E$ is the support of some free $\ol{T}$-controlled $R$-module.
\end{prop}

In the proof of this proposition we shall use the following

\begin{lem}\label{estonoacabanunca}
Given an injective controlled homomorphism $\varphi\colon
R\grupo{A}_\alpha\rightarrow R\grupo{B}_\beta$ we have that
$\alpha(A)'\subset\beta(B)'$.
\end{lem}

\begin{proof}
For any $e\in \alpha(A)'$ we can take a sequence
$\set{a_n}_{n\in\N}\subset A$ with
$\lim\limits_{n\rightarrow\infty}\alpha(a_n)=e$. Since $\varphi$
is injective $\operatorname{carr}(\varphi(a_n))$ is non-empty for
every $n\in\N$ so we can take elements $b_n\in
\operatorname{carr}(\varphi(a_n))$. By definition of controlled
homomorphism $\lim\limits_{n\rightarrow\infty}\beta(b_n)=e$, hence
$e\in\beta(B)'$ and the inclusion holds.
\end{proof}

\begin{proof}[Proof of (\ref{clasiiso})]
The case $R=\Z$ and $\ol{T}$ a tree-like space follows from
\cite{iht} III.4.8 and III.4.16. In general condition (1) is
necessary since an isomorphism of free $\ol{T}$-controlled
$R$-modules is also an isomorphism between the underlying
$R$-modules. Moreover, condition (2) is necessary by
(\ref{estonoacabanunca}). By Remark \ref{zzz} it is enough to make
the proof for tree-like spaces. In the rest of the proof we shall
suppose that $\ol{T}$ is tree-like.

If $\alpha(A)'=\alpha(B)'=\emptyset$ then $A$ and $B$ are both
finite and any isomorphism $R\grupo{A}\simeq R\grupo{B}$ is a
controlled isomorphism $R\grupo{A}_\alpha\simeq R\grupo{B}_\beta$.
If $\alpha(A)'=\alpha(B)'\neq\emptyset$ then $A$ and $B$ are
infinite countable, so $\Z\grupo{A}\simeq\Z\grupo{B}$. Since the
proposition holds for $R=\Z$ there is a controlled isomorphism
$\varphi\colon\Z\grupo{A}_\alpha\simeq\Z\grupo{B}_\beta$. Now one
can check that $\varphi\otimes R\colon R\grupo{A}_\alpha\simeq
R\grupo{B}_\beta$ is an isomorphism of free $\ol{T}$-controlled
$R$-modules.

Finally, given $K\subset E$ compact, if $\Z\grupo{C}_\gamma$ is a
free $\ol{T}$-controlled $\Z$-module whose support is $K$ then the
support of $R\grupo{C}_\gamma$ is $K$ as well, because it only
depends on $\gamma$. The proof is now complete.
\end{proof}

The next result follows directly from Remark \ref{zzz},
Proposition \ref{clasiiso} and the definition of controlled
homomorphisms.

\begin{prop}\label{desc}
Up to isomorphism, the endomorphism algebra of a free
$\ol{T}$-controlled $R$-module with support $E$ only depends on
$E$. Moreover, it is isomorphic to $R(E)$.
\end{prop}

The last isomorphism of this proposition is given by the fact that
the basis of a free $\ol{T}$-controlled $R$-module
$R\grupo{A}_\alpha$ with support $E$ must be infinite countable,
and hence it can be identified with the non-negative integers
$A=\NN$. Moreover, we can suppose that $\alpha$ is the inclusion
of a discrete subspace $\alpha\colon A\subset Y$, changing $A$ by
$\alpha(A)$ if necessary. Now we are in the same situation as in
the beginning of the introduction. We also derive from
(\ref{desc}) that the isomorphism class of the $R$-algebra $R(E)$
only depends on $E$, as we claimed in the introduction.

\begin{prop}
Every free $\ol{T}$-controlled $R$-module is a retract of any
object whose support is $E$.
\end{prop}

\begin{proof}
Recall from (\ref{clasiiso}) that all objects with support $E$ are
isomorphic. By (\ref{zzz}) it is enough to check the proposition
for $\ol{T}$ a tree-like space. For $R=\Z$ and $\ol{T}$ tree-like
this proposition is contained in the proof of \cite{iht} V.3.4.
The result for arbitrary rings follows from the special case
$R=\Z$. More precisely, given a free $\ol{T}$-controlled
$R$-module $R\grupo{A}_\alpha$ if the support of
$R\grupo{B}_\beta$ is $E$ then $\Z\grupo{A}_\alpha$ is a retract
of $\Z\grupo{B}_\beta$ (the supports only depend on the height
functions) hence we obtain a retraction of $R\grupo{B}_\beta$ onto
$R\grupo{A}_\alpha$ by tensoring by $R$.
\end{proof}





As a consequence of this proposition we get by (\ref{mitchell})
the following equivalence of categories which will be used from
now on as an identification.

\begin{cor}\label{unaeq}
The evaluation functor in a free $\ol{T}$-controlled $R$-module
with support $E$ induces an additive equivalence of abelian
categories $\mo(\M_R(\ol{T}))\simeq\mo(R(E))$ which restricts to
another one $\fp({\M_R(\ol{T})})\simeq \fp(R(E))$.
\end{cor}

\begin{rem}\label{mismo}
If $R\simeq R^{op}$, in particular if $R$ is commutative, the
transposition of matrices and an explicit isomorphism $R\simeq
R^{op}$ induce isomorphisms of $R$-ringoids
$\M_R(\ol{T})\simeq\M_R(\ol{T})^{op}$ (preserving objects) and
$R$-algebras $R(E)\simeq R(E)^{op}$, compare \cite{canc}, so in
this case right-modules over $\M_R(\ol{T})$ or $R(E)$ are the same
as left-modules.
\end{rem}

In the following proposition we compute the dimension of the
$k$-algebra $k(E)$, $k$ any field.

\begin{prop}\label{dim}
$\dim k(E)=2^{\aleph_0}$.
\end{prop}

\begin{proof}
The $k$-vector space of all $\NN\times \NN$ matrices is the direct
product of $\NN\times \NN$ copies of $k$, and it is known that
$\dim\Pi_{\NN\times \NN}k=2^{\aleph_0}$, hence $\dim k(E)\leq
2^{\aleph_0}$. Moreover, for any element $(\mat{a}_n)_{n\in
\NN}\in\prod_{\NN} k$ the diagonal matrix
$(\mat{b}_{ij})_{i,j\in\NN}$ with $\mat{b}_{nn}=\mat{a}_n$ belongs
to $k(E)$, therefore $k(E)$ has a vector subspace isomorphic to
$\prod_{\NN} k$, and $\dim\prod_{\NN} k=2^{\aleph_0}$ as well, so
the equality of the statement holds.
\end{proof}

\subsection{The special case $\card E$ finite}\label{finite}
If $\card{E}$ is finite, since $E$ is metrizable, it must have the
discrete topology, so $E$ is the discrete set $\n{n}$ with
$n=\card E$ points. For this particular space we can take a
tree-like space $\ol{T}_n=(\hat{T}_n,T_n,\n{n})$ where $T_n$ is a
locally compact tree with $n$ Freudenthal ends and $\hat{T}_n$ is
the Freudenthal compactification of $T_n$, see \cite{iht} III.1.3.
Moreover, if $T^0_n$ is the vertex set of $T_n$ and $\delta\colon
T^0_n\subset T_n$ is the inclusion, the support of
$R\grupo{T^0_n}_\delta$ is $\n{n}$, in particular $R(\n{n})$ is
the endomorphism ring of this object. Let us fix the following
particular tree $T_n$: the vertex set of $T_n$ is
$$T_n^0=\set{v_0}\cup\set{v_m^1,\dots v_m^n}_{m\geq 1},$$ and there
are edges joining $v_0$ with $v_1^i$ and $v_m^i$ with $v^i_{m+1}$
$(1\leq i\leq n, m\geq 1)$. The additive category $\M_R(\ol{T}_n)$
is equivalent to the full subcategory of objects
$R\grupo{A}_\alpha$ such that $\alpha(A)\subset T^0_n$, compare
the proof of \cite{iht} V.3.4. From now on we shall always work in
this subcategory, and we denote it by $\M_R(\ol{T}_n)$ as well.

We are going to give an alternative description for controlled
homomorphisms in $\M_R(\ol{T}_n)$. For this we define the
following sets for any height function $\alpha\colon A\rightarrow
T^0_n\subset T_n$ $(1\leq i\leq n, j\geq 1)$
$$A^i_j=\bigcup_{l\geq j}\alpha^{-1}(v^i_l)$$
A morphism $\varphi\colon
 R\grupo{B}_\beta\rightarrow R\grupo{A}_\alpha$ in $\M_R(\ol{T}_n)$ is controlled
if and only if for every $m\geq 1$ there exists $M\geq 1$ such
that $\varphi(B^i_M)\subset R\grupo{A^i_m}$ for any $1\leq i\leq
n$. We shall omit the superindex $i$ when $n=1$. Moreover, for the
next sections we fix the following notation $(m\geq 0)$
$${}_mA=\alpha^{-1}(v_0)\cup\left[\bigcup_{1\leq i\leq n}\bigcup_{l\leq m}\alpha^{-1}(v^i_l)\right],\;\; {}_mA^i_j={}_mA\cap A^i_j.$$

\begin{rem}\label{matricillas}
If $n=1$, $T_1=[0,+\infty)$ is the half-line and $T^0_1=\NN$ the
non-negative integers. Moreover $R(\n{1})$ is the $R$-algebra
$\mathrm{RCFM}(R)$ of $\NN\times\NN$ matrices with entries in $R$
such that every row and every column has a finite number of
non-zero entries (\emph{row-column-finite matrices}), compare
\cite{iht} V.3.8.
\end{rem}


\begin{rem}\label{reduc}
If $E$ is any compact metrizable space with at least $n$ points we
can fully include $\M_R(\ol{T}_n)$ into $\M_R(\ol{T})$. For this
we only need to take $n$ disjoint sequences $\set{v^i_m}_{m\geq
1}$ $(1\leq i\leq n)$ contained in $Y$ converging in $X$ to $n$
different points belonging to $E$, and an additional point $v_0\in
Y$ out of the sequences. Now we identify $\M_R(\ol{T}_n)$ with the
full subcategory of $\M_R(\ol{T})$ given by objects
$R\grupo{A}_\alpha$ with $\alpha(A)\subset
\set{v_0}\cup\set{v_m^1,\dots v_m^n}_{m\geq 1}$. If we call
$\uf{F}$ to this full inclusion, we get another one
$\uf{F}_*\colon\mo(\M_R(\ol{T}_n))\rightarrow\mo(\M_R(\ol{T}))$
together with a retraction up to natural equivalence
$\uf{F}^*\colon\mo(\M_R(\ol{T}))\rightarrow\mo(\M_R(\ol{T}_n))$.
Moreover, the first functor $\uf{F}_*$ restricts to the full
subcategories of f. p. modules, see Section \ref{rym}, hence the
decomposition problem for f. p. $R({\n{n}})$-modules is included
in the decomposition problem for f. p. $R(E)$-modules, in
particular we only need to prove Theorem \ref{rep} for $\card E$
finite.
\end{rem}

\section{Countably presented $R$-modules as finitely presented
$\mathrm{RCFM}(R)$-modules}\label{fijate}\label{count}

There is a full exact inclusion of abelian categories
$$\mathfrak{i}\colon\mo(R)\To\mo(\mathrm{RCFM}(R))$$ defined by
$\mathfrak{i}{M}=\hom_R(R\grupo{\NN},{M})$. The
ring $\mathrm{RCFM}(R)$ acts on $\mathfrak{i}{M}$ by
endomorphisms of $R\grupo{\NN}$.

Let $\mathfrak{f}\colon\M_R(\ol{T}_1)\rightarrow\mo(R)$ be the
forgetful functor  which sends a free $\ol{T}_1$-controlled
$R$-module to its underlying $R$-module. The
$\mathrm{RCFM}(R)$-module $\mathfrak{i}{M}$ can be
regarded as the functor
$\mathfrak{i}{M}=\hom_R(\mathfrak{f},{M})\colon\M_R(\ol{T}_1)^{op}\rightarrow\Ab$.

\begin{prop}\label{adj}
The functor $\mathfrak{i}$ has an exact left-adjoint
$\mathfrak{r}$ such that $\mathfrak{ri}$ is naturally equivalent
to the identity functor. Moreover, $\mathfrak{r}$ can be chosen to
be the evaluation functor in a free $\ol{T}_1$-controlled
$R$-module with one generator.
\end{prop}

\begin{proof}
Let $R\grupo{e}_\phi$ be a free $\ol{T}_1$-controlled $R$-module
whose basis is a singleton $\set{e}$ (all these objects are
isomorphic in $\M_R(\ol{T}_1)$ by (\ref{clasiiso})). The
endomorphism ring of this object is $R$, hence the evaluation
functor $ev_{R\grupo{e}_\phi}$ takes values in the category of
$R$-modules. Let us see that the exact functor
$ev_{R\grupo{e}_\phi}$ is left-adjoint to $\mathfrak{i}$. A
left-adjoint for $\mathfrak{i}$ exists and can be constructed by
left additive Kan extension, see \cite{rso} 6, so we have just to
check that
$\hom_{\M_R(\ol{T}_1)}(R\grupo{A}_\alpha,\mathfrak{i}{M})=\hom_R(ev_{R\grupo{e}_\phi}R\grupo{A}_\alpha,{M})$
for any free $\ol{T}_1$-controlled $R$-module $R\grupo{A}_\alpha$
in a natural way. This follows from the obvious natural
identification $ev_{R\grupo{e}_\phi}R\grupo{A}_\alpha=R\grupo{A}$
and Yoneda's lemma.
\end{proof}

\begin{cor}\label{in}
If ${M}$ is an $R$-module and $\mathcal{N}$ an
$\mathrm{RCFM}(R)$-module then there are natural isomorphisms
$(n\geq 0)$
$$\ext^n_{R(E)}(\mathcal{N},\mathfrak{i}{M})\simeq\ext^n_R(\mathfrak{r}\mathcal{N},{M}).$$
In particular if $R=k$ is a field the $\mathrm{RCFM}(k)$-modules
$\mathfrak{i}{M}$ are all injective.
\end{cor}

An $R$-module is \emph{countably presented} provided it is the
cokernel of a morphism between free $R$-modules with countable
basis. Obviously the cokernel of a morphism between countably
presented $R$-modules is countably presented as well. In
particular direct summands of countably presented $R$-modules are
countably presented.

\begin{prop}\label{igno}
The functor $\mathfrak{i}$ sends countably presented $R$-modules
to finitely presented $\mathrm{RCFM}(R)$-modules.
\end{prop}

In the proof of this proposition we shall use the
row-column-finite matrices $\mat{A}$ and $\mat{B}$ defined by
\begin{itemize}
\item $\mat{a}_{i+1,i}=1$ $(i\in\NN)$ and $\mat{a}_{ij}=0$ in
other cases,

\item $\mat{b}_{\frac{n(n+1)}{2}+i,\frac{(n-1)n}{2}+i}=1$ for any
$n>0$ and $0\leq i<n$, and $\mat{b}_{ij}=0$ otherwise.
\end{itemize}
And we regard $\mathrm{RCFM}(R)$ as the endomorphism $R$-algebra
of the free $\ol{T}_1$-controlled $R$-module $R\grupo{\NN}_\delta$
where $\delta\colon\NN\subset[0,+\infty)$ is the inclusion, see
Remark \ref{matricillas} and Proposition \ref{desc}.

\begin{proof}[Proof of (\ref{igno})]
Since $\mathfrak{i}$ is exact it will be enough to check the
proposition for the countably presented $R$-modules $R$ and
$R\grupo{\NN}$. Recall that
$\hom_{\M_R(\ol{T}_1)}(R\grupo{\NN}_\delta,\mathfrak{i}{M})=
\hom_R(R\grupo{\NN},{M})$ for any $R$-module
${M}$. The $\mathrm{RCFM}(R)$-modules $\mathfrak{i}R$ and
$\mathfrak{i}R\grupo{\NN}$ are the cokernels of
$(\mat{I}-\mat{A})$ and $(\mat{I}-\mat{B})$ respectively. The
natural projections onto the cokernel are given by the
homomorphisms $R\grupo{\NN}\rightarrow R$ and
$R\grupo{\NN}\rightarrow R\grupo{\NN}$ defined on generators by
$n\mapsto 1$ $(n\geq 0)$ and $\frac{n(n+1)}{2}+i\mapsto i$ $(n\geq
i\geq 0)$ respectively.
\end{proof}

As a consequence of Propositions \ref{adj} and \ref{igno}, and
Remark \ref{reduc} we get the next

\begin{cor}
The representation problem for countably generated $R$-modules is
contained into the representation problem for finitely presented
$R(E)$-modules, where $E$ is any non-empty compact metrizable
space.
\end{cor}

\section{A review on pro-categories}\label{repro}

Any partially ordered set (poset) $\Lambda$ can be regarded as a
small category with a unique morphism $u\rightarrow v$ provided
$u\geq v$, $u,v\in\Lambda$. A poset $\Lambda$ is \emph{directed}
if given $u,v\in \Lambda$ there exists $w\in \Lambda$ with $w\geq
u,v$. Moreover, $\Lambda$ is \emph{cofinite} if the set $\set{u\in
\Lambda\,;\,u\leq v}$ is finite for every $v\in \Lambda$.

A \emph{pro-object} or \emph{inverse system} $X_\bullet$ over a
category $\mathbf{C}$ is a functor $X_\bullet\colon
\Lambda\rightarrow\mathbf{C}$ where $\Lambda$ is a directed
cofinite poset. If $u\in\Lambda$ we usually write
$X_u=X_\bullet(u)$. The morphisms $X_\bullet (u\rightarrow v)$
$(u,v\in \Lambda, u\geq v)$ are the \emph{bonding morphisms} of
$X_\bullet$, and $\Lambda$ is the \emph{indexing set} of the
inverse system.

The category $\mathrm{pro}-\mathbf{C}$ has objects inverse systems
over $\mathbf{C}$. Morphism sets are given by the following
formula
\begin{equation}\label{promor}
\hom_{\mathrm{pro}-\mathbf{C}}(X_\bullet,Y_\bullet)=
\lim\limits_{v}\colim\limits_{u} \hom_{\mathbf{C}}(X_u,Y_v).
\end{equation}

We identify any object in $\mathbf{C}$ with the inverse system
whose indexing set is a singleton $\Lambda=*$. This defines a full
inclusion of categories $\mathbf{C}\subset
\mathrm{pro}-\mathbf{C}$. This inclusion has a right-adjoint, the
(inverse) limit functor
$\lim\colon\mathrm{pro}-\mathbf{C}\rightarrow\mathbf{C}$, $\lim
X_\bullet=\lim_{u}X_u$.

The category $\mathrm{pro}-\mathbf{C}$ is abelian whenever
$\mathbf{C}$ is, see \cite{cp} 6.4. This will be always the case,
because we are only going to use in this context the category
$\mathbf{C}=\mo(k)$ of $k$-vector spaces. If $V$ is a vector space
and $X_\bullet$ an inverse system of vector spaces, then by
(\ref{promor})
\begin{equation*}
\hom_{\mathrm{pro}-\mo(k)}(V,X_\bullet)=\lim\limits_v
\hom_k(V,X_v)= \hom_k(V,\lim  X_\bullet).
\end{equation*}
Hence, since $\hom_k(V,-)$ is an exact functor in the category of
vector spaces, the Grothendieck spectral sequence (see
\cite{hilton} 9.3) yields an isomorpshim
\begin{equation}\label{gro}
\ext^1_{\mathrm{pro}-\mo(k)}(V,X_\bullet)=\hom_k(V,{\lim }^1
X_\bullet).
\end{equation}

\section{Numerical invariants of finitely presented
$k(\n{n})$-modules}\label{nume}

In this section we shall define invariants of the isomorphism
class of a f. p. $k(\n{n})$-module lying in the abelian monoids
$\N_{\infty,n}$ $(n\geq 1)$. The abelian monoid $\N_{\infty,n}$
has $n+1$ generators
$$1,\;\infty_1,\;\dots\;\infty_n$$ and $2n$ relations
$$1+\infty_i=\infty_i,\;\infty_i+\infty_i=\infty_i,\;\;\;(1\leq i\leq
n).$$ As a set $\N_{\infty,n}$ is
$$\N_{\infty,n}=\NN\sqcup\set{\infty_S=\sum_{i\in S}\infty_i\,;\,\emptyset\neq
S\subset\set{1,\dots n}}.$$ For the sake of simplicity if $n=1$ we
write $\N_{\infty,1}=\N_\infty$ and $\infty_1=\infty$. 

Let $\varphi\colon k\grupo{B}_\beta\rightarrow k\grupo{A}_\alpha$
be a morphism in $\M_k(\ol{T}_n)$. We define the element
$\lambda_\varphi\in\N_{\infty,n}$ in the following way: if the
next vector space is finite-dimensional
\begin{equation*}
L_\varphi=\frac{k\grupo{A}}{\bigcap_{m\geq 1}\sum_{i=1}^n
\left[k\grupo{A_m^i}+\varphi(k\grupo{B})\right]},
\end{equation*} then $\lambda_\varphi=\dim L_\varphi$,
otherwise $\lambda_\varphi=\infty_S$, where $S\subset\set{1,\dots
n}$ is the biggest subset such that if $i\notin S$ then there
exists $M\geq 1$ with $k\grupo{A^i_M}\subset\varphi(k\grupo{B})$.

\begin{prop}
The element $\lambda_\varphi$ only depends on the isomorphism
class of the f. p. $k(\n{n})$-module $\coker\varphi$, and
$\lambda_{\varphi\oplus\psi}=\lambda_\varphi+\lambda_\psi$.
\end{prop}

\begin{proof}
One can readily check by using the alternative description of
controlled homomorphisms given in Subsection \ref{finite} that the
correspondence $\varphi\mapsto \Upsilon(\varphi)=L_\varphi$
defines an additive functor $\Upsilon$ from
$\pair(\M_k(\ol{T}_n))$ to the category of $k$-vector spaces.
Moreover, if $V^{\varphi,i}_\bullet$ is the inverse system of
$k$-vector spaces indexed by $\N$ and given by $(1\leq i\leq n)$
$$V^{\varphi,i}_m=\frac{k\grupo{A^i_m}+\varphi(k\grupo{B})}{\varphi(k\grupo{B})},$$
and the obvious inclusions as bonding morphisms, the
correspondences
$\varphi\mapsto\Theta_i(\varphi)=V^{\varphi,i}_\bullet$ also
define additive functors $\Theta_i$ from $\pair(\M_k(\ol{T}_n))$
to the pro-category of pro-vector spaces. Furthermore, it is easy
to see that the functors $\Upsilon$, $\Theta_i$ $(1\leq i\leq n)$
factor through the natural equivalence relation $\sim$ in
$\pair(\M_k(\ol{T}_n))$, hence the first statement of the
proposition follows from (\ref{cokereq}) and the fact that
$\lambda_\varphi$ is defined as $\dim \Upsilon(\varphi)$ provided
this vector space is finite-dimensional, and otherwise
$\lambda_\varphi=\infty_S$ where $S=\set{i\in\set{1,\dots
n}\,;\,\Theta_i(\varphi)\simeq\!\!\!\!\!\!/\;\; 0}$. The second
part of the statement follows from the additivity of the functors
$\Upsilon$, $\Theta_i$ $(1\leq i\leq n)$.
\end{proof}

If the following vector space has finite dimension,
\begin{equation*}
M_\varphi^i=\frac{\bigcap_{m\geq 1}\left[k\grupo{A^i_m}+
\varphi(k\grupo{B})\right]}{\bigcap_{m\geq
1}\left\{\left[k\grupo{A^i_m}+
\varphi(k\grupo{B})\right]\cap\left[\sum_{j\neq i} k\grupo{A_m^j}+
\varphi(k\grupo{B})\right]\right\}},
\end{equation*} the element
$\mu_\varphi^i\in\N_{\infty}$ $(1\leq i\leq n)$ is defined as
$\mu_\varphi^i=\dim M_\varphi^i$, otherwise
$\mu^i_\varphi=\infty$.

\begin{prop}
The elements $\mu^i_\varphi$ $(1\leq i\leq n)$ only depend on the
isomorphism class of the f. p. $k(\n{n})$-module $\coker\varphi$,
and $\mu_{\varphi\oplus\psi}^i=\mu^i_\varphi+\mu^i_\psi$.
\end{prop}

\begin{proof}
By using the characterization of controlled homomorphisms given in
Subsection \ref{finite} one readily checks that the
correspondences $\varphi\mapsto M_\varphi^i$ define additive
functors from $\pair(\M_k(\ol{T}_n))$ to the category of
$k$-vector spaces. Moreover, these functors factor through the
natural equivalence relation $\sim$, hence the proposition follows
from (\ref{cokereq}).
\end{proof}

In order to define elements $\nu^i_\varphi\in\N_{\infty}$ $(1\leq
i\leq n)$ we introduce inverse systems of $k$-vector spaces
$U^{\varphi,i}_\bullet$, indexed by the set $\N\times\N$ with the
product partial order, given by
\begin{equation*}
U^{\varphi,i}_{pq}=\frac{k\grupo{A^i_p}\cap\varphi(k\grupo{B})}
{k\grupo{A^i_p}\cap\varphi(k\grupo{B^i_q})},
\end{equation*} and bonding homomorphisms induced by the obvious inclusions of vector spaces.
If the limit of $U^{\varphi,i}_\bullet$ is finite-dimensional we
set $\nu^i_\varphi=\dim\lim  U^{\varphi,i}_\bullet,$ otherwise
$\nu^i_\varphi=\infty$.

\begin{prop}
The elements $\nu^i_\varphi$ $(1\leq i\leq n)$ only depend on the
isomorphism class of the f. p. $k(\n{n})$-module $\coker\varphi$,
and $\nu^i_{\varphi\oplus\psi}=\nu^i_\varphi+\nu^i_\psi$.
\end{prop}

\begin{proof}
One can check by using the description of controlled homomorphisms
given in Subsection \ref{finite} that the correspondences
$\varphi\mapsto U^{\varphi,i}_\bullet$ are additive functors from
$\pair(\M_k(\ol{T}_n))$ to the pro-category of pro-vector spaces,
and these functors factor through the natural equivalence relation
$\sim$, hence the proposition follows from (\ref{cokereq}).
\end{proof}

The former propositions of this section are summarized in the
following

\begin{cor}\label{resumen}
There are well defined morphisms of abelian monoids $(n\in\N)$
\begin{equation*}
\Phi_n\colon\operatorname{Iso}(\fp(k(\n{n})))\To\N_{\infty,n}\times\prod_{i=1}^n\N_{\infty}
\times\prod_{i=1}^n\N_{\infty}
\end{equation*} which send the isomorphism class $[\um{M}]$ of a f. p. $k(\n{n})$-module
$\um{M}=\coker\varphi$ to
\begin{equation*}
\Phi_n([\um{M}])=\left(\lambda_\varphi,\left(\mu_\varphi^i\right)_{i=1}^n,
\left(\nu_\varphi^i\right)_{i=1}^n\right).
\end{equation*}
\end{cor}

From now on we shall write $\lambda_\um{M}=\lambda_\varphi$,
$\mu^i_\um{M}=\mu^i_\varphi$ and $\nu^i_\um{M}=\nu^i_\varphi$
$(1\leq i\leq n)$ if $\um{M}=\coker\varphi$ and omit the
superindex $i$ when $n=1$.

\begin{rem}\label{teincluyo}
There are $n$ full inclusions
$\uf{F}^i\colon\mathbf{M}_k(\ol{T}_1)\rightarrow\mathbf{M}_k(\ol{T}_n)$
$(1\leq i\leq n)$ defined by identifying $T^0_1=\NN$ (see
(\ref{matricillas})) with the subset
$\set{v_0}\cup\set{v^i_m}_{m\geq 1}\subset T^0_n$ in the obvious
way, see (\ref{reduc}).
\end{rem}

The next proposition can be easily checked by using the
commutativity of (\ref{conmuta}) and the right-exactness of the
functors $\uf{F}^i_*$.

\begin{prop}\label{metiendo}
If $\um{M}$ is a f. p. $k(\n{1})$-module then for every $1\leq
i\leq n$
\begin{itemize}
\item $\lambda_{\uf{F}^i_*\um{M}}=\lambda_\um{M}$ if
$\lambda_\um{M}\in\NN$, and $\lambda_{\uf{F}^i_*\um{M}}=\infty_i$
otherwise,

\item $\mu^i_{\uf{F}^i_*\um{M}}=\mu_\um{M}$ and
$\mu^j_{\uf{F}^i_*\um{M}}=0$ if $j\neq i$,

\item $\nu^i_{\uf{F}^i_*\um{M}}=\nu_\um{M}$ and
$\nu^j_{\uf{F}^i_*\um{M}}=0$ if $j\neq i$.
\end{itemize}
\end{prop}

\section{Classification of finitely presented
$\mathrm{RCFM}(k)$-modules}\label{rcfm}

Recall from (\ref{matricillas}) that the $k$-algebra $k(\n{1})$
coincides with $\mathrm{RCFM}(k)$, the $k$-algebra of matrices
$\mat{A}=(\mat{a}_{ij})_{i,j\in\NN}$ with entries in $k$ such that
every row and every column has at most a finite number of
non-trivial entries (row-column-finite matrices). Those matrices
are the endomorphisms of the free $\ol{T}_1$-controlled $k$-vector
space $k\grupo{\NN}_\delta$, where
$\delta\colon{\NN\subset[0,+\infty)}$ is the inclusion of the
vertex set. The unit element of the $k$-algebra $\mathrm{RCFM}(k)$
is the identity matrix $\mat{I}$ with $\mat{i}_{ii}=1$ $(i\in\NN)$
and $\mat{i}_{ij}=0$ if $i\neq j$. For the sake of simplicity we
abbreviate $\um{R}=\mathrm{RCFM}(k)$.

Consider the matrices $\mat{A}$ and $\mat{B}$ used in the proof of
(\ref{igno}), they are defined as
\begin{itemize}
\item $\mat{a}_{i+1,i}=1$ $(i\in\NN)$ and $\mat{a}_{ij}=0$ in
other cases,

\item $\mat{b}_{\frac{n(n+1)}{2}+i,\frac{(n-1)n}{2}+i}=1$ for any
$n>0$ and $0\leq i<n$, and $\mat{b}_{ij}=0$ otherwise.
\end{itemize}
We define the following $\um{R}$-modules
\begin{itemize}
\item $\um{A}=\frac{\um{R}}{\mat{A} \um{R}}$,

\item $\um{B}=\frac{\um{R}}{(\mat{I}-\mat{A}) \um{R}}$,

\item $\um{C}=\frac{\um{R}}{(\mat{I}-\mat{A}^t) \um{R}}$,

\item $\um{B}_\infty=\frac{\um{R}}{(\mat{I}-\mat{B}) \um{R}}$,

\item $\um{C}_\infty=\frac{\um{R}}{(\mat{I}-\mat{B}^t) \um{R}}$.
\end{itemize}

The main result of this section is the next

\begin{thm}[Classification of f. p. $\mathrm{RCFM}(k)$-modules]\label{clasifR}
There is a solution to the decomposition problem in the category
of f. p. $\mathrm{RCFM}(k)$-modules given by the following
elementary modules
 $$\mathcal{A},\; \um{R},\;
\mathcal{B},\; \mathcal{B}_\infty,\;\mathcal{C},\;
\mathcal{C}_\infty,$$ and elementary isomorphisms
$$\mathcal{A}\oplus \um{R}\simeq \um{R},\;\,
\um{R}\oplus \um{R}\simeq \um{R},\;\, \mathcal{B}\oplus
\mathcal{B}_\infty\simeq \mathcal{B}_\infty,\;\,
$$ $$ \mathcal{B}_\infty\oplus\mathcal{B}_\infty\simeq
\mathcal{B}_\infty,\;\, \mathcal{C}\oplus \mathcal{C}_\infty\simeq
\mathcal{C}_\infty,\;\, \mathcal{C}_\infty\oplus
\mathcal{C}_\infty\simeq \mathcal{C}_\infty.$$

\end{thm}

This theorem implies Theorems \ref{rep} and \ref{generel} for
$\card E=1$. It is a direct consequence of the next two results.
We shall use the following notation, given $n\in\NN$ we write
$\um{A}_n$, $\um{B}_n$ and $\um{C}_n$ for the direct sum of $n$
copies of $\um{A}$, $\um{B}$ or $\um{C}$, respectively, and
$\um{A}_\infty=\um{R}$.

\begin{thm}\label{iny}%
For every f. p. $\um{R}$-module $\um{M}$ there is an isomorphism
$$\um{M}\simeq\um{A}_{\lambda_\um{M}}\oplus\um{B}_{\mu_\um{M}}\oplus\um{C}_{\nu_\um{M}}.$$
\end{thm}

\begin{prop}\label{imagenes}
We have the following equalities:
\begin{itemize}
\item $\Phi_1(\um{A})=(1,0,0)$,

\item $\Phi_1(\um{B})=(0,1,0)$,

\item $\Phi_1(\um{C})=(0,0,1)$,

\item $\Phi_1(\um{R})=(\infty,0,0)$,

\item $\Phi_1(\um{B}_\infty)=(0,\infty,0)$,

\item $\Phi_1(\um{C}_\infty)=(0,0,\infty)$.
\end{itemize}
\end{prop}

In particular we have that

\begin{cor}\label{esiso}
The monoid morphism
$$\Phi_1\colon\operatorname{Iso}(\fp(k(\n{1})))\To\N_{\infty}\times\N_{\infty}
\times\N_{\infty}$$ is an isomorphism.
\end{cor}

The proof of Proposition \ref{imagenes} will be given later.
Theorem \ref{iny} is a direct consequence of Lemmas \ref{alfa},
\ref{yvan2}, \ref{casielfinal} and \ref{es0}. These are some of
the hardest technical results of this paper. In fact the rest of
this long section is highly technical. It is focused towards
proving Theorem \ref{iny}, although some interesting corollaries
on the homological algebra of finitely presented $\um{R}$-modules
are derived from the technical lemmas. These homological results
are used in the proof of (\ref{iny}) as well as in the appendix.
We advice the reader to skip this material in a first reading.

\medskip

The next result is an easy computation.

\begin{lem}\label{mult}
The following equalities hold in $\um{R}$:
\begin{enumerate}
\item $\mat{A}^t\mat{A}=\mat{I}$,

\item $\mat{B}^t\mat{B}=\mat{I}$.
\end{enumerate}
\end{lem}

The next lemma follows directly from (\ref{clasiiso}).

\begin{lem}\label{unaclas}
Two free $\ol{T}_1$-controlled $k$-vector spaces
$k\grupo{A}_\alpha$, $k\grupo{B}_\beta$ are isomorphic in
$\M_k(\ol{T}_1)$ if and only if $A$ and $B$ have the same
cardinal.
\end{lem}


\begin{lem}\label{ya}
The $\um{R}$-module $\um{A}$ is isomorphic to any 1-dimensional
free $\ol{T}_1$-controlled $k$-vector space.
\end{lem}

\begin{proof}
If $k\grupo{e}_\phi$ is a 1-dimensional free $\ol{T}_1$-controlled
$k$-vector space, the cokernel of $\mat{A}$ is given by the
controlled homomorphism $\varphi\colon
k\grupo{\NN}_\delta\twoheadrightarrow k\grupo{e}_\phi$ defined
over the basic elements as $0\mapsto e$ and $n\mapsto 0$ for
$n>0$.
\end{proof}

The proof of Proposition \ref{imagenes} is as follows.

\begin{proof}[Proof of \ref{imagenes}]
In this proof we omit some straightforward but tedious
computations which can be carried out by the interested reader
with not too much difficulty. We shall write $\N_p$ $(p\geq 1)$
for the set of naturals $\geq p$.

The $\um{R}$-module $\um{R}$ corresponds to the free
$\ol{T}_1$-controlled $k$-vector space $k\grupo{\NN}_\delta$ where
$\delta\colon\NN\subset[0,+\infty)$ is the inclusion, hence it is
the cokernel of the trivial morphism $0\colon 0\rightarrow
k\grupo{\NN}_\delta$, and the equalities $\lambda_\um{R}=\infty$,
$\mu_\um{R}=0$ hold immediately, moreover,
$U^0_{m,n}=k\grupo{\N_m}$ for all $m,n\in\N$, and given $M\geq m$,
$N\geq n$ the corresponding bonding homomorphism in $U^0_\bullet$
is the inclusion $U^0_{M,N}=k\grupo{\N_M}\subset k\grupo{\N_m}=
U^0_{m,n}$, therefore $\lim
U^0_\bullet=\bigcap\limits_{m\in\N}k\grupo{\N_m}=0$ and
$\nu_\um{R}=0$.

By (\ref{ya}) the $\M_k(\ol{T}_1)$-module $\um{A}$ is the cokernel
of the trivial morphism $0\rightarrow k\grupo{e}_\phi$ where
$k\grupo{e}_\phi$ is a 1-dimensional free $\ol{T}_1$-controlled
$k$-vector space, hence the equality $\Phi_1(\um{A})=(1,0,0)$
follows easily.

One can check that
$k\grupo{\N_n}+(\mat{I}-\mat{A})(k\grupo{\NN})=k\grupo{\NN}$ for
all $n\in\N$, and
$k\grupo{\NN}/(\mat{I}-\mat{A})(k\grupo{\NN})\simeq k$ generated
by the class of any $n\in\NN$, so $\lambda_{\um{B}}=0$ and
$\mu_{\um{B}}=1$. Moreover,
$k\grupo{\N_n}\cap(\mat{I}-\mat{A})(k\grupo{\NN})=(\mat{I}-\mat{A})(k\grupo{\N_n})$
and hence $U^{(\mat{I}-\mat{A})}_{n,n}=0$ for all $n\in\N$,
therefore $\lim U^{(\mat{I}-\mat{A})}_\bullet=0$ since the
diagonal subset $\set{(n,n)\,;\,n\in\N}\subset\N\times\N$ is
cofinal, so $\nu_{\um{B}}=0$.

One can see that $(\mat{I}-\mat{A}^t)(k\grupo{\NN})=k\grupo{\NN}$,
hence $\lambda_{\um{C}}=0=\mu_{\um{C}}$, moreover
$(\mat{I}-\mat{A}^t)(k\grupo{\N_n})$ is generated by the set
$\set{m-(m-1)}_{m\geq n}$, therefore
$U^{(\mat{I}-\mat{A}^t)}_{n-1,n}\simeq k$ generated by the class
of any $m\geq n-1$ and the bonding homomorphism
$U^{(\mat{I}-\mat{A}^t)}_{n,n+1}\rightarrow
U^{(\mat{I}-\mat{A}^t)}_{n-1,n}$ is an isomorphism, so again by
cofinality we see that $\lim U^{(\mat{I}-\mat{A}^t)}_\bullet\simeq
k$, in particular $\nu_{\um{C}}=1$.

The vector space $k\grupo{\N_n}+(\mat{I}-\mat{B})(k\grupo{\NN})$
is the whole $k\grupo{\NN}$, so $\lambda_{\um{B}_\infty}=0$,
moreover, a basis of
$k\grupo{\NN}/(\mat{I}-\mat{B})(k\grupo{\NN})$ is
$\set{\frac{n(n+3)}{2}}_{n\in\NN}$, hence
$\mu_{\um{B}_\infty}=\infty$. One can check that
$k\grupo{\N_n}\cap(\mat{I}-\mat{B})(k\grupo{\NN})=(\mat{I}-\mat{B})(k\grupo{\N_n})$,
therefore $U^{(\mat{I}-\mat{B})}_{n,n}=0$, $\lim
U^{(\mat{I}-\mat{B})}_\bullet=0$ and $\nu_{\um{B}_\infty}=0$.

Finally $(\mat{I}-\mat{B}^t)(k\grupo{\NN})=k\grupo{\NN}$, so
$\lambda_{\um{C}_\infty}=0=\mu_{\um{C}_\infty}$, and there are
isomorphisms $(n\in\NN)$
\begin{equation*}
U^{(\mat{I}-\mat{B}^t)}_{\frac{(n+1)(n+2)}{2},\frac{n(n+1)}{2}}\simeq
k\left\langle\frac{n(n+1)}{2},\,\dots\,
\frac{n(n+1)}{2}+n\right\rangle,
\end{equation*} moreover, the following bonding homomorphism $(n>0)$
\begin{equation*}
U^{(\mat{I}-\mat{B}^t)}_{\frac{(n+1)(n+2)}{2},\frac{n(n+1)}{2}}\To
U^{(\mat{I}-\mat{B}^t)}_{\frac{n(n+1)}{2},\frac{(n-1)n}{2}}
\end{equation*} sends $\frac{n(n+1)}{2}+m$ to $\frac{(n-1)n}{2}+m$
if $m<n$ and $\frac{n(n+1)}{2}+n$ to the trivial element, so $\lim
U^{(\mat{I}-\mat{B}^t)}_\bullet=\prod_{\NN}k$ is the direct
product of an infinite countable number of copies of $k$ and the
equality $\nu_{\um{C}_\infty}=\infty$ holds.
\end{proof}

\begin{lem}\label{otroiso}
There is an $\um{R}$-module isomorphism $\frac{\um{R}}{\mat{B}
\um{R}}\simeq\um{R}$.
\end{lem}

\begin{proof}
Let $A\subset\NN$ be the infinite subset
$A=\set{\frac{n(n+3)}{2}}_{n\in\NN}$, and $\alpha\colon
A\subset\NN$ the inclusion. The next sequence, where $\varphi$ is
the obvious projection, is exact
$$k\grupo{\NN}_\delta\st{\mat{B}}\hookrightarrow k\grupo{\NN}_\delta
\st{\varphi}\twoheadrightarrow k\grupo{A}_\alpha.$$ Hence the
lemma follows from (\ref{unaclas}).
\end{proof}

\begin{lem}\label{semeten}
Left-multiplication by one of the following matrices induces an
injective right-$\um{R}$-module homomorphism
$\um{R}\rightarrow\um{R}$,
$$\mat{A},(\mat{I}-\mat{A}),(\mat{I}-\mat{A}^t),(\mat{I}-\mat{B}),(\mat{I}-\mat{B}^t).$$
\end{lem}

\begin{proof}
The matrix $\mat{A}$ has a left-inverse in $\um{R}$ by
(\ref{mult}). The other matrices have a left-inverse either in the
$k$-algebra $\mathrm{CFM}(k)$ of column-finite matrices or in the
$k$-algebra $\mathrm{RFM}(k)$ of row-finite matrices. Both
$k$-algebras contain $\um{R}$, moreover,
$\um{R}=\mathrm{CFM}(k)\cap\mathrm{RFM}(k)$. The $k$-algebra
$\mathrm{CFM}(k)$ is just the endomorphism ring of the $k$-vector
space $k\grupo{\NN}$, and there is an isomorphism
$\mathrm{RFM}(k)\simeq\mathrm{CFM}(k)^{op}$ given by
transposition. More precisely, let $\mat{C}$, $\mat{D}$ be the
matrices in $\mathrm{RFM}(k)$ defined by $\mat{c}^1_{ij}=1$ if
$i\geq j$ and zero otherwise, and
$\mat{d}_{\frac{(m-1)m}{2}+i,\frac{(n-1)n}{2}+i}=1$ $(m\geq
n>i\geq 0)$ and trivial in other cases. One can check that
$\mat{C}(\mat{I}-\mat{A})=\mat{I}$,
$\mat{C}^t(\mat{I}-\mat{A}^t)=\mat{I}$,
$\mat{D}(\mat{I}-\mat{B})=\mat{I}$ and
$\mat{D}^t(\mat{I}-\mat{B}^t)=\mat{I}$, hence the lemma follows.
\end{proof}

\begin{prop}\label{short}
There are extensions of $\um{R}$-modules
\begin{enumerate}
\item $\um{A}\hookrightarrow\um{B}\twoheadrightarrow\um{C}$, \item
$\um{R}\hookrightarrow\um{B}_\infty\twoheadrightarrow\um{C}_\infty$.
\end{enumerate}
\end{prop}

\begin{proof}
By using (\ref{mult}), (\ref{otroiso}) and (\ref{semeten}) we get
the following equalities, isomorphisms and short exact sequences,
which correspond to the extensions of the statement
$$\frac{\um{R}}{\mat{A} \um{R}}\simeq
\frac{(\mat{I}-\mat{A}^t) \um{R}}{(\mat{I}-\mat{A}^t) \mat{A}
\um{R}}= \frac{(\mat{I}-\mat{A}^t) \um{R}}{(\mat{I}-\mat{A})
\um{R}}\hookrightarrow\frac{\um{R}} {\mat{(\mat{I}-\mat{A})}
\um{R}}\twoheadrightarrow\frac{\um{R}} {\mat{(\mat{I}-\mat{A}^t)}
\um{R}},$$
$$\um{R}\simeq\frac{\um{R}}{\mat{B} \um{R}}\simeq
\frac{(\mat{I}-\mat{B}^t) \um{R}}{(\mat{I}-\mat{B}^t) \mat{B}
\um{R}}= \frac{(\mat{I}-\mat{B}^t) \um{R}}{(\mat{I}-\mat{B})
\um{R}}\hookrightarrow\frac{\um{R}} {\mat{(\mat{I}-\mat{B})}
\um{R}}\twoheadrightarrow\frac{\um{R}} {\mat{(\mat{I}-\mat{B}^t)}
\um{R}}.$$
\end{proof}

The proof of the following proposition is contained in the proof
of (\ref{igno}).

\begin{prop}\label{isoin}
Given a $k$-vector space $V$, if $\dim V<\aleph_0$ then
$\mathfrak{i}V=\um{B}_{\dim V}$, and $\mathfrak{i}V=\um{B}_\infty$
if $\dim V=\aleph_0$.
\end{prop}


The next corollary follows from (\ref{isoin}) and (\ref{in}).

\begin{cor}\label{medioin}
The $\um{R}$-module $\um{B}_d$ is injective for every
$d\in\N_\infty$.
\end{cor}

In the next lemma we show that one can adapt the basis of a
countably generated vector space to a decreasing filtration.

\begin{lem}\label{noname}
Let $V_0\supset V_1\supset\cdots\supset V_n\supset
V_{n+1}\supset\cdots$ be a decreasing sequence of $k$-vector
spaces such that $V_0$ is the union of an increasing sequence of
finite-dimensional subspaces $V^0_0\subset
V^1_0\subset\cdots\subset V^n_0\subset V_0^{n+1}\subset\cdots$,
$V_0=\bigcup_{n\in\NN} V_0^n$. If we set $V^{-1}_0=0$,
$V^n_m=V^n_0\cap V_m$ $(n+1,m\in\NN)$,
$V_\infty=\bigcap_{n\in\NN}V_n$ and choose (finite and possibly
empty) sets $\set{a^l_{n m}\,;\,1\leq l\leq r_{n m}}\subset V^n_m$
such that the sets
$$\set{a^l_{n m}+(V^{n-1}_m+V^n_{m+1})\,;\,1\leq l\leq r_{nm}}$$
are basis of $V^n_m/(V^{n-1}_m+V^n_{m+1})$ $(n,m\in\NN)$, then
given $m,n,p\in\NN$ with $m\leq p$
\begin{enumerate}
\item $\set{a^{l_{ij}}_{ij}+V^n_p\,;\, i\leq n, m\leq j< p, 1\leq
l_{ij}\leq r_{ij}}$ is a basis of $V^n_m/V^n_p$,

\item $\set{a^{l_{ij}}_{ij}+V_\infty\,;\, i\leq n, m\leq j, 1\leq
l_{ij}\leq r_{ij}}$ is a basis of $(V^n_m+V_\infty)/V_\infty$,

\item $\set{a^{l_{ij}}_{ij}+V_\infty\,;\, i\in\NN, m\leq j, 1\leq
l_{ij}\leq r_{ij}}$ is a basis of $V_m/V_\infty$.
\end{enumerate}
\end{lem}

\begin{proof}
Since $V^n_0$ is a finite-dimensional vector space it is artinian
and the decreasing sequence $V^n_0\supset
V^n_1\supset\cdots\supset V^n_m\supset V^n_{m+1}\supset\cdots$
stabilizes, it is, there exists $M_n\in\NN$ such that
$V^n_m=V^n_{M_n}$ for every $m\geq M_n$, in particular
$V^n_{M_n}=V^n_0\cap V_\infty$.
If we choose for every $n\in\NN$ the minimum $M_n$ satisfying this
condition then $M_n\leq M_{n+1}$ since
$$V^n_{M_{n+1}}=V^n_{M_{n+1}}\cap V^{n+1}_{M_{n+1}}=V^n_0\cap V_{M_{n+1}}\cap V^{n+1}_0\cap V_\infty=
V^n_0\cap V_\infty=V^n_{M_n}.$$

Notice that $(1)$ is trivial for $m\geq M_n$ since
$V^n_m=V^n_p=V^n_{M_n}$ and $\set{a_{ij}^l\,;\,1\leq l\leq
r_{ij}}=\emptyset$ whenever $p\geq m\geq M_n$, $i\leq n$ and
$j\geq M_n$. Therefore the unique elements $(n,m,p)$ for which we
still have to check $(1)$ lie in the set $S=\set{(n,m,p)\,;\,
n\in\NN,0\leq m\leq M_n,p\geq m}$. Let us order this set in the
following way
$$(n,m,p)\leq(n',m',p')\Leftrightarrow\left\{
\begin{array}{ll}
    n<n' & \\
     \hbox{or}&  \\
     n=n' \;\hbox{and} \;m>m'& \\
      \hbox{or} &\\
     n=n',\,m=m'\;\hbox{and}\;p\leq p'.&
\end{array}%
\right.    $$ One readily checks that this is a well order on $S$,
since the second coordinate has an upper bound (depending on the
first one). The minimum of $S$ is $(0,M_0,M_0)$, moreover if
$m<M_n$ the element $(n,m,m)$ is the least upper bound of
$\set{(n,m+1,p)\,;\,p>m}$, and given $n>0$ the element
$(n,M_n,M_n)$ is the least upper bound of the set
$\set{(n-1,m,p)\,;\,m\leq M_{n-1},p\geq m}$. Any other element in
$S$ is a successor. We have already checked $(1)$ for the elements
$(n,M_n,p)\in S$, moreover, it is trivial for $(n,m,m)\in S$,
hence $(1)$ holds for the minimum and all limit elements in $S$. A
generic successor in $S$ has the form $(n,m,p+1)$ for some $p\geq
m$. Notice that we have already check $(1)$ for some successors as
well, namely for those with $m=M_n$. We are now going to proceed
by induction, it is, we shall prove $(1)$ for every successor in
$S$ with $m<M_n$ supposing that $(1)$ holds for all the strictly
lower elements. We are going to distinguish three cases:

\underline{For $(0,m,p+1)$} $(1)$ follows from the exactness of
the sequence
$$\frac{V^0_p}{V^0_{p+1}}\hookrightarrow \frac{V^0_m}{V^0_{p+1}}\twoheadrightarrow \frac{V^0_m}{V^0_p},$$
the equality $V^{-1}_0=0$ and the inequality $(0,m,p)<(0,m,p+1)$,
and the inequality $(0,p+1,p)<(0,m,p+1)$ if $p<M_0$ or the
equalities $V_p^0=V^0_{p+1}=V^0_{M_0}$ and
$\set{a_{0p}^l\,;\,1\leq l\leq r_{0p}}=\emptyset$ if $p\geq M_0$.


\underline{For $(n,m,m+1)$ with $n>0$
} $(1)$ is a consequence of the exactness of the sequence
$$\frac{V^{n-1}_m+V^n_{m+1}}{V^n_{m+1}}\hookrightarrow \frac{V^n_m}{V^n_{m+1}}
\twoheadrightarrow \frac{V^n_m}{V^{n-1}_m+V^n_{m+1}},$$ the
obvious isomorphism
$$\frac{V^{n-1}_m}{V^{n-1}_{m+1}}\simeq\frac{V^{n-1}_m+V^n_{m+1}}{V^n_{m+1}},$$
the inequality $(n-1,m,m+1)<(n,m,m+1)$ if $m<M_{n-1}$, or the
equalities $V^{n-1}_m=V^{n-1}_{m+1}=V^{n-1}_{M_{n-1}}$ and
$\set{a_{in}\,;\, 1\leq l\leq r_{in}}=\emptyset$ if $i\leq n-1$
and $m\geq M_{n-1}$.

\underline{For $(n,m,p+1)$ with $n>0$ and $p>m$} $(1)$ follows
from the exactness of the sequence
$$\frac{V^n_p}{V^n_{p+1}}\hookrightarrow \frac{V^n_m}{V^n_{p+1}}\twoheadrightarrow
\frac{V^n_m}{V^n_p},$$ the inequality $(n,m,p)<(n,m,p+1)$, and the
inequality $(n,p,p+1)<(n,m,p+1)$ if $p< M_n$ or the equalities
$V^n_p=V^n_{p+1}=V^n_{M_n}$ and $\set{a_{ip}^l\,;\,1\leq l\leq
r_{ip}}=\emptyset$ if $i\leq n$ and $p\geq M_n$.

Once we have seen that $(1)$ holds, $(2)$ is a consequence of
$(1)$ for $p=M_n$, the isomorphism
$(V^n_m+V_\infty)/V_\infty\simeq V^n_m/V^n_{M_n}$, and the fact
that $\set{a_{ij}^l\,;\,1\leq l\leq r_{ij}}=\emptyset$ is the
empty set for $i\leq n$ and $j\geq M_n$. Finally $(3)$ follows
from $(2)$ and the equality $V_m=\bigcup_{n\in\NN}V^n_m$.
\end{proof}

The next proposition is an interesting consequence of the previous
lemma. It does not hold in general when the ground ring is not a
field, compare \cite{projdim}.

\begin{prop}\label{loslibres}
The image of a morphism between finitely generated free
$\M_k(\ol{T}_1)$-modules is finitely generated free.
\end{prop}

\begin{proof}
Let $\varphi\colon k\grupo{B}_\beta\rightarrow k\grupo{A}_\alpha$
be a morphism in $\mathbf{M}_k(\ol{T}_1)$. If we define the
$k$-vector spaces $V_0=\varphi(k\grupo{B})$,
$V_n=\varphi(k\grupo{B_n})\subset k\grupo{A}$ $(n\in\N)$ and
$V_0^n=\varphi(k\grupo{{}_nB})$ $(n\geq 0)$ we can apply Lemma
\ref{noname}. Moreover, with the notation of that lemma
$V_\infty=0$ since for every $n\geq 1$ there exists $N_n\geq 1$
such that $V_{N_n}\subset k\grupo{A_n}$ and $\cap_{n\geq
1}k\grupo{A_n}=0$. We define the set
$\ul{B}=\set{a_{ij}^l\,;\,i,j\in \NN,1\leq l\leq r_{ij}}$ and the
function $\ul{\beta}\colon \ul{B}\rightarrow\NN\subset
[0,+\infty)$ by $\ul{\beta}(a_{nm}^l)=m$. By (\ref{noname}) (3)
the set $\ul{B}$ is a basis of $V_0$ and $\ul{B}_m$ a basis of
$V_m$ $(m\geq 1)$ since $V_\infty=0$. The function $\beta$ is a
height function because the cardinal of $\ul{\beta}^{-1}(m)$ is
$\dim V_m/V_{m+1}<\aleph_0$ $(m\in\NN)$. Moreover, the inclusion
$k\grupo{\ul{B}}=\varphi(k\grupo{B})\subset k\grupo{A}$ and the
projection
$k\grupo{B}\twoheadrightarrow\varphi(k\grupo{B})=k\grupo{\ul{B}}$
give rise to controlled homomorphisms
$k\grupo{\ul{B}}_{\ul{\beta}}\hookrightarrow k\grupo{A}_\alpha$
and $k\grupo{B}_\beta\twoheadrightarrow
k\grupo{\ul{B}}_{\ul{\beta}}$ which are an $\M_k(\ol{T}_1)$-module
monomorphism and epimorphism respectively and their composition is
$\varphi$, hence $k\grupo{\ul{B}}_{\ul{\beta}}$ together with
these morphisms is the image of $\varphi$.
\end{proof}



\begin{cor}\label{presentacion}
Any finitely presented $\M_k(\ol{T}_1)$-module is the cokernel of
a monomorphism between finitely generated free
$\M_k(\ol{T}_1)$-modules.
\end{cor}

\begin{cor}\label{pd}
Finitely presented $\um{R}$-modules have projective dimension
$\leq 1$.
\end{cor}

\begin{cor}\label{medioin2}
We have $\ext^1(\um{M},\um{C}_d)=0$ for any f. p. $\um{R}$-module
$\um{M}$ and $d\in\N_\infty$.
\end{cor}

\begin{proof}
By (\ref{pd}) the functor $\ext^1(\um{M},-)$ is right-exact, hence
the corollary follows from (\ref{short}) and (\ref{medioin}).
\end{proof}

Now we begin with the lemmas which prove Theorem \ref{iny}.

\begin{lem}\label{alfa}
Given any f. p. $\um{R}$-module $\um{M}$, there exists another f.
p. $\um{R}$-module $\um{N}$ with $\lambda_\um{N}=0$ such that
$\um{M}\simeq\um{A}_{\lambda_\um{M}}\oplus\um{N}$.
\end{lem}

\begin{proof}\renewcommand{\theequation}{\alph{equation}}\setcounter{equation}{0}
Suppose that $\um{M}=\coker\varphi$ for some $\varphi\colon
k\grupo{B}_\beta\rightarrow k\grupo{A}_\alpha$ in
$\M_k(\ol{T}_1)$. Let us consider the decreasing sequence of
$k$-vector spaces given by $V_0=k\grupo{A}/\varphi(k\grupo{B})$
and
$$V_n=\frac{k\grupo{A_n}+\varphi(k\grupo{B})}{\varphi(k\grupo{B})},\,\,\, n\in\N.$$
The vector space $V_0$ is the union of the following sequence of
finite-dimensional $k$-vector spaces $(n\in\NN)$
$$V^n_0=\frac{k\grupo{{}_nA}+\varphi(k\grupo{B})}{\varphi(k\grupo{B})}.$$

If $\set{a^l_{n m}\,;\,1\leq l\leq r_{n m}}\subset V^n_m$ is a set
as in (\ref{noname}) we can suppose that
$a_{nm}^l=e_{nm}^l+\varphi(k\grupo{B})$ for some $e_{nm}^l\in
k\grupo{{}_nA_m}$, here we use the next obvious isomorphism
$$\frac{k\grupo{{}_nA_m}}{k\grupo{{}_nA_m}\cap \varphi(k\grupo{B})}\simeq
\frac{k\grupo{{}_nA_m}+\varphi(k\grupo{B})}{\varphi(k\grupo{B})}=V^n_m.$$

We consider the set $C=\set{a_{nm}^{l_{nm}}+V_\infty\,;\,
n,m\in\NN,1\leq l_{nm}\leq r_{nm}}$ and the function $\gamma\colon
C\rightarrow\NN\subset[0,+\infty)$ with
$\gamma(a_{nm}^{l_{nm}}+V_\infty)=m$. This function is a height
function, since the set
$\gamma^{-1}(m)=\set{a_{nm}^{l_{nm}}+V_\infty\,;\,n\in\NN,1\leq
l_{nm}\leq r_{nm}}$ is bijective with a basis of $V_m/V_{m+1}$ by
(\ref{noname}) (3), and we have the following surjection and
isomorphisms
$$\frac{k\grupo{A_m}}{k\grupo{A_{m+1}}}\twoheadrightarrow
\frac{k\grupo{A_m}}{k\grupo{A_{m+1}}+\left[\varphi(k\grupo{B})\cap
k\grupo{A_m}\right]}
\simeq\frac{k\grupo{A_m}+\varphi(k\grupo{B})}{k\grupo{A_{m+1}}+\varphi(k\grupo{B})}\simeq
\frac{V_m}{V_{m+1}},$$ and $\dim
k\grupo{A_m}/k\grupo{A_{m+1}}=\card\alpha^{-1}(m)<\aleph_0$. The
underlying $k$-vector space of $k\grupo{C}_{\gamma}$ is
$V_0/V_\infty$, moreover, the natural projection
\begin{equation}\label{raca}
k\grupo{A}\twoheadrightarrow\frac{k\grupo{A}}{\bigcap\limits_{n\geq
1}\left[k\grupo{A_n}
+\varphi(k\grupo{B})\right]}\simeq\frac{V_0}{V_\infty}=k\grupo{C}
\end{equation}
give rise to a $\ol{T}_1$-controlled homomorphism $\upsilon_0
\colon k\grupo{A}_\alpha\rightarrow k\grupo{C}_{\gamma}$ with
$\upsilon_0 \varphi=0$, hence $\upsilon_0$ induces a morphism
$\upsilon\colon \um{M}\rightarrow k\grupo{C}_{\gamma}$.
Furthermore, the section $V_0/V_\infty\hookrightarrow k\grupo{A}$
which sends $a^l_{nm}+V_\infty$ to $e^l_{nm}$ determines another
$\ol{T}_1$-controlled homomorphism $\tau_0\colon
k\grupo{C}_{\gamma}\rightarrow k\grupo{A}_\alpha$ with $\upsilon_0
\tau_0 =1$, in particular if $\tau\colon
k\grupo{C}_{\gamma}\rightarrow \um{M}$  is the morphism induced by
$\tau_0$ we have that $\upsilon\tau=1$, hence $\um{M}\simeq
k\grupo{C}_{\gamma}\oplus \um{N}$ where $\um{N}=\coker\tau$.
Notice that the morphism $(\varphi,\tau_0)\colon
k\grupo{B}_\beta\oplus k\grupo{C}_{\gamma}\rightarrow
k\grupo{A}_\alpha$ is a finite presentation of $\um{N}$, and by
(\ref{raca}) we have the following equality and inclusions for
every $m\geq 1$
$$k\grupo{A}=\bigcap\limits_{n\geq 1}\left[k\grupo{A_n}
+\varphi(k\grupo{B})\right]\oplus \tau (k\grupo{C})\subset
k\grupo{A_m}+\varphi(k\grupo{B})+\tau (k\grupo{C})\subset
k\grupo{A},$$ therefore $\lambda_\um{N}=0$.

Observe that by (\ref{unaclas}) and (\ref{ya})
$k\grupo{C}_{\gamma}$ is isomorphic to the free
$\ol{T}_1$-controlled $k$-vector space which corresponds to
$\um{R}$ provided $\lambda_\um{M}=\infty$, and to the direct sum
of $\lambda_\um{M}$ copies of $\um{A}$ otherwise.
\end{proof}

\begin{lem}\label{yvan2}
Given a f. p. $\um{R}$-module $\um{M}$ with $\lambda_\um{M}=0$,
there exists another f. p. $\um{R}$-module $\um{N}$ with
$\lambda_\um{N}=\mu_\um{N}=0$ such that
$\um{M}\simeq\um{B}_{\mu_\um{M}}\oplus\um{N}$.
\end{lem}

\begin{proof}\renewcommand{\theequation}{\alph{equation}}\setcounter{equation}{0}
Suppose that $\um{M}$ is the cokernel of $\varphi\colon
k\grupo{B}_\beta\rightarrow k\grupo{A}_\alpha$ in
$\M_k(\ol{T}_1)$. Since $\lambda_\um{M}=0$ we have that
$k\grupo{A}=k\grupo{A_m}+\varphi(k\grupo{B})$ for every $m\in\N$.
Let $V_\bullet$ be the inverse system indexed by $\N\times\N$
given by
\begin{equation*}
V_{mn}=\frac{k\grupo{A_m}}{k\grupo{A_m}\cap
\varphi(k\grupo{B_n})},
\end{equation*} and bonding homomorphisms induced by the obvious
inclusions of vector spaces. There are inclusions
$U_{mn}^\varphi\subset V_{mn}$ with quotients
\begin{equation*}
\frac{k\grupo{A_m}}{k\grupo{A_m}\cap \varphi(k\grupo{B})}\simeq
\frac{k\grupo{A_m}+\varphi(k\grupo{B})}{\varphi(k\grupo{B})}=
\frac{k\grupo{A}}{\varphi(k\grupo{B})}.
\end{equation*} This determines a short exact sequence in the
pro-category of pro-vector spaces
\begin{equation}\label{unacorta}
U^\varphi_\bullet\hookrightarrow
V_\bullet\twoheadrightarrow\frac{k\grupo{A}}{\varphi(k\grupo{B})}.
\end{equation}
Here we regard ${k\grupo{A}}/{\varphi(k\grupo{B})}$ as the inverse
system indexed by a singleton.

The vector space $U^\varphi_{mn}$ is always finite dimensional,
because it is contained in
$\varphi(k\grupo{B})/\varphi(k\grupo{B_n})\simeq\varphi(k\grupo{{}_{n-1}B})$
and ${}_{n-1}B$ is a finite set. Since finite-dimensional vector
spaces are artinian it is easy to see that $U^\varphi_\bullet$
satisfies the Mittag-Leffler property, in particular $\lim^1
U^\varphi_\bullet=0$ and hence by (\ref{gro})
$\ext^1({k\grupo{A}}/{\varphi(k\grupo{B})}, U^\varphi_\bullet)=0$,
so the sequence (\ref{unacorta}) admits a splitting $s\colon
{k\grupo{A}}/{\varphi(k\grupo{B})}\hookrightarrow V_\bullet$. This
splitting is given by splittings $s_{mn}\colon
{k\grupo{A}}/{\varphi(k\grupo{B})}\hookrightarrow V_{mn}$ of the
natural projections
$V_{mn}\twoheadrightarrow{k\grupo{A}}/{\varphi(k\grupo{B})}$ which
are compatible with the bonding homomorphisms of $V_\bullet$.

Let $\tilde{C}$ be a basis of $k\grupo{A}/\varphi(k\grupo{B})$.
This basis is either finite
$\tilde{C}=\set{b_1,\,\dots\,b_{\mu_\um{M}}}$ if
$\mu_\um{M}\in\NN$, or infinite countable
$\tilde{C}=\set{b_n}_{n\in\NN}$ if $\mu_\um{M}=\infty$. Moreover,
since $\varphi$ is controlled there exists an increasing sequence
of natural numbers $\set{l_n}_{n\geq 1}$ with
$\varphi(k\grupo{B_{l_n}})\subset k\grupo{A_n}$. We choose
elements $b^{n-1}_m\in k\grupo{A_n}$ and $y^{n-1}_m\in
k\grupo{B_{l_n}}$ such that
$b^{n-1}_m+\varphi(k\grupo{B_{l_n}})=s_{n,l_n}(b_m)\in
V_{n,l_n}=k\grupo{A_n}/\varphi(k\grupo{B_{l_n}})$ and
$\varphi(y^{n-1}_m)=b^{n}_m-b^{n-1}_m$ for every $n\in\N$ and $m$
in the corresponding range. Furthermore, we define the sets
${}^nC\subset k\grupo{A_{n+1}}$ and $C$ in the following way:
${}^nC=\{b^n_1,\,\dots\,b^n_{\mu_\um{M}}\}$ and
$C=\coprod_{n\in\NN}{}^nC$ if $\mu_\um{M}\in\NN$, and
${}^nC=\set{b^n_0,\,\dots\,b^n_n}\cup\set{b^m_m\,;\,m>n}$ and
$C=\bigcup_{n\in\NN}{}^nC=\set{b^n_m\,;\,n\geq m\geq 0}$ if
$\mu_\um{M}=\infty$. Let $\gamma\colon C\rightarrow\NN\subset
[0,+\infty)$ be the height function given by $\gamma(b^n_m)=n$ and
$\psi$ the endomorphism of $k\grupo{C}_\gamma$ given by
$\psi(b^n_m)=b^{n+1}_m-b^n_m$.

One readily checks that
$\coker\psi=\mathfrak{i}(k\grupo{A}/\varphi(k\grupo{B}))$ and the
natural projection
$k\grupo{C}_\gamma\twoheadrightarrow\mathfrak{i}(k\grupo{A}/\varphi(k\grupo{B}))$
is given by the homomorphism $p_1\colon k\grupo{C}\rightarrow
k\grupo{A}/\varphi(k\grupo{B})=k\grupo{\tilde{C}}$ defined by
$p_1(b^n_m)=b_m$. For this one uses the finite presentations
constructed in the proof of (\ref{igno}) and, if
$\mu_{\mathcal{M}}=\infty$, the bijection $\NN\approx C$ which
sends $m\in\NN$, with $\frac{n(n-1)}{2}\leq m<\frac{(n+1)n}{2}$
for some $n\in\NN$, to $b^{n-1}_{m-\frac{n(n-1)}{2}}$. Moreover,
by (\ref{isoin}) $\coker\psi=\mathcal{B}_{\mu_\mathcal{M}}$.

The homomorphism $\tau_0\colon k\grupo{C}\rightarrow k\grupo{A}$
induced by the inclusions ${}^nC\subset k\grupo{A_{n+1}}\subset
k\grupo{A}$ determines a controlled homomorphism $\tau_0\colon
k\grupo{C}_\gamma\rightarrow k\grupo{A}_\alpha$. Moreover, the
homomorphism $\tau_1\colon k\grupo{C}\rightarrow k\grupo{B}$ given
by $\tau_1(b^n_m)=y^n_m$ defines a controlled homomorphism
$\tau_1\colon k\grupo{C}_\gamma\rightarrow k\grupo{B}_\beta$ with
$\varphi\tau_1=\tau_0\psi$, hence $\tau_0$ gives rise to a
$\M_k(\ol{T}_1)$-module morphism
$\tau\colon\um{B}_{\mu_\um{M}}\rightarrow \mathcal{M}$.

Let us check that $\tau$ is a monomorphism. Given a free
$\ol{T}_1$-controlled vector space $k\grupo{D}_\phi$ Yoneda's
lemma yields a natural identification
$\hom_\mathcal{R}(k\grupo{D}_\phi,\mathfrak{i}(k\grupo{A}/\varphi(k\grupo{B})))=\hom_k(k\grupo{D},k\grupo{A}/\varphi(k\grupo{B}))$.
This identification carries a morphism $\upsilon\colon
k\grupo{D}_\phi\rightarrow\mathfrak{i}(k\grupo{A}/\varphi(k\grupo{B}))$
represented by $\upsilon_0\colon k\grupo{D}_\phi\rightarrow
k\grupo{C}_\gamma$ to the vector space homomorphism
$p_1\upsilon_0$. If $p_2\colon k\grupo{A}\twoheadrightarrow
k\grupo{A}/\varphi(k\grupo{B})$ is the natural projection then
$p_2\varphi=0$ and $p_1=p_2\tau_0$. Moreover $\tau\upsilon=0$ if
and only if  $\tau_0\upsilon_0=\varphi\eta$ for some controlled
homomorphism $\eta\colon k\grupo{D}_\phi\rightarrow
k\grupo{B}_\beta$, so in this case
$p_1\upsilon_0=p_2\tau_0\upsilon_0=p_2\varphi\eta=0$, it is,
$\upsilon=0$, therefore $\tau$ is a monomorphism.

By (\ref{medioin}) if $\um{N}=\coker\tau$ then
$\um{M}\simeq\um{B}_{\mu_\um{M}}\oplus\um{N}$. The morphism
$(\varphi,\tau_0)\colon k\grupo{B}_\beta\oplus
k\grupo{C}_\gamma\rightarrow k\grupo{A}_\alpha$ is a finite
presentation of $\um{N}$ and by construction
$\varphi(k\grupo{B})+\tau_0(k\grupo{C})=k\grupo{A}$, it is
$\lambda_\um{N}=\mu_\um{N}=0$.
\end{proof}

\begin{lem}\label{casielfinal}
Let $\um{M}$ be a f. p. $\um{R}$-module with
$\lambda_\um{M}=\mu_\um{M}=0$ then there exists another one
$\um{N}$ with $\lambda_\um{N}=\mu_\um{N}=\nu_\um{N}=0$ such that
$\um{M}\simeq\um{C}_{\nu_\um{M}}\oplus\um{N}$.
\end{lem}

In the proof of this lemma we shall use the next one.

\begin{lem}\label{uniso}
Let $\set{d_n}_{n\in\N}$ be an increasing sequence of integers
with $\lim\limits_{n\rightarrow\infty}d_n=\infty$. Consider the
set $A=\set{(n,m)\,;\,n\in\N,m\leq d_n}\subset\N\times\N$, the
height  function $\alpha\colon A\rightarrow\NN\subset[0,+\infty)$
with $\alpha(n,m)=n$, and the endomorphism $\varphi$ of
$k\grupo{A}_\alpha$ with $\varphi(n,m)=(n,m)$ if $n=1$ or $n>1$
and $m>d_{n-1}$, and $\varphi(n,m)=(n,m)-(n-1,m)$ otherwise. Then
$\coker\varphi$ is isomorphic to $\um{C}_\infty$.
\end{lem}

\begin{proof}
Consider the infinite countable subsets $$A_1=\set{(1,m)\,;\,1\leq
m\leq d_1}\cup\set{(n,m)\,;\,n>1, d_{n-1}<m\leq d_n}\subset A,$$
$A_2=A-A_1$, $B_1=\set{(n+1)n/2}_{n\in\NN}\subset\NN$ and
$B_2=\NN-B_1$. The lexicographic order from the left on $A$ is a
well order without limit elements, since the second coordinate of
an element $(n,m)\in A$ is bounded by $d_n$, hence the restriction
of this order to the subsets $A_1$ and $A_2$ induces enumerations
$A_1=\set{e_1^n}_{n\in\NN}$, $A_2=\set{e^2_n}_{n\in\NN}$.
Similarly the usual order in $\NN$ induces enumerations in the
subsets $B_1=\set{f_1^n}_{n\in\NN}$ and
$B_2=\set{f_2^n}_{n\in\NN}$. Now the theorem follows from the
bijection $\NN\approx A$ which sends $f^n_i$ to $e^n_i$
$(i=1,2;n\in\NN)$.
\end{proof}

\begin{proof}[Proof of (\ref{casielfinal})]\renewcommand{\theequation}{\alph{equation}}\setcounter{equation}{0}
If $\um{M}=\coker[\varphi\colon k\grupo{B}_\beta\rightarrow
k\grupo{A}_\alpha]$ the equalities $\lambda_\um{M}=\mu_\um{M}=0$
are equivalent to $\varphi(k\grupo{B})=k\grupo{A}$. Let
$\phi\colon \lim U^\varphi_\bullet \rightarrow U^\varphi_\bullet$
be the canonical pro-morphism. This pro-morphism is given by
vector space homomorphisms $\phi_{mn}\colon \lim U^\varphi_\bullet
\rightarrow U^\varphi_{mn}$ compatible with the bonding
homomorphisms of $U^\varphi_\bullet$. Since $\varphi$ is
controlled there is an increasing sequence of natural numbers
$\set{m_n}_{n\geq 1}$ such that $\varphi(k\grupo{B_{m_n}})\subset
k\grupo{A_n}$.

If $\nu_\um{M}\in\NN$  and $\set{a_1,\dots a_{\nu_\um{M}}}$ is a
basis of $\lim U^\varphi_\bullet$ we define
${}^nC=\set{a^n_1,\dots a^n_{\nu_\um{M}}}\subset k\grupo{A_n}$ as
a set such that
$\phi_{n,m_n}(a_i)=a^n_i+\varphi(k\grupo{B_{m_n}})$ $(1\leq
i\leq\nu_\um{M})$, and choose elements $y^n_i\in
k\grupo{B_{m_{n-1}}}$ if $n>1$ and $y^1_i\in k\grupo{B}$ with
$a^n_i-a^{n-1}_i=\varphi(y^n_i)$ $(n>1)$ and
$a^1_i=\varphi(y^1_i)$. If $\nu_\um{M}=\infty$ we take
${}^nC=\set{a^n_i}_{i=1}^{d_n}\subset k\grupo{A_n}$  such that
$\set{a^n_i+\varphi(k\grupo{B_{m_n}})}$ is a basis of
$\phi_{n,m_n}(\lim U^\varphi_\bullet)$, here we use that
$U^\varphi_\bullet$ is an inverse system of finite-dimensional
vector spaces, compare the proof of (\ref{yvan2}). The bonding
homomorphisms of $U^\varphi_\bullet$ induce surjections
$\phi_{n+1,m_{n+1}}(\lim
U^\varphi_\bullet)\twoheadrightarrow\phi_{n,m_n}(\lim
U^\varphi_\bullet)$, hence $d_n\leq d_{n+1}$ and we can suppose
without loose of generality that there exist $y^n_i\in
k\grupo{B_{m_{n-1}}}$ $(n>1)$ and $y^1_i\in k\grupo{B}$ such that
$a^n_i-a^{n-1}_i=\varphi(y^n_i)$ if $n>1$ and $i\leq d_{n-1}$, and
$a^n_i=\varphi(y^n_i)$ if $n>1$ and $d_{n-1}< i\leq d_n$ or $n=1$
and $i\leq d_1$.

We define the height function $\gamma\colon C=\coprod_{n\geq
1}{}^nC\rightarrow\NN\subset [0,+\infty)$ as $\gamma(a^n_i)=n$,
and the controlled homomorphisms $\tau_0\colon
k\grupo{C}_\gamma\rightarrow k\grupo{A}_\alpha$, $\tau_1\colon
k\grupo{C}_\gamma\rightarrow k\grupo{B}_\beta$, $\psi\colon
k\grupo{C}_\gamma\rightarrow k\grupo{C}_\gamma$ by
$\tau_0(a_i^n)=a^n_i$, $\tau_1(a_i^n)=y^n_i$, and
$\psi(a^n_i)=a^n_i-a^{n-1}_i$ if $n>1$, and $\nu_\um{M}\in\NN$ or
$\nu_\um{M}=\infty$ and $i\leq d_{n-1}$, and $\psi(a^n_i)=a^n_i$
otherwise. One can readily check, by using the bijection
$\N\approx\NN\colon n\mapsto n-1$ if $\nu_\um{M}\in\NN$ or
(\ref{uniso}) if $\nu_\um{M}=\infty$, that $\coker\psi\simeq
\um{C}_{\nu_\um{M}}$. Moreover,
$(\tau_1,\tau_0)\colon\psi\rightarrow\varphi$ is a morphism in
$\pair(\M_k(\ol{T}_1))$ which induces an $\um{R}$-module morphism
$\tau\colon \um{C}_{\nu_\um{M}}\rightarrow\um{M}$.

In order to check that $\tau$ is a monomorphism of
$\M_k(\ol{T}_1)$-modules we are going to prove that $\phi$ is a
monomorphism of pro-vector spaces if $\nu_\um{M}\in\NN$. In this
case $\ker\phi$ is an inverse system of finite-dimensional vector
spaces, in particular it satisfies the Mittag-Leffler property. If
we apply the left-exact functor $\lim$ to the exact sequence of
pro-vector spaces
$$\ker\phi\hookrightarrow\lim U^\varphi_\bullet\st{\phi}\To
U^\varphi_\bullet$$ we get another one
$$\lim\ker\phi\hookrightarrow\lim\lim
U^\varphi_\bullet\st{=}\To \lim U^\varphi_\bullet,$$ so
$\lim\ker\phi=0$ and hence $\ker\phi=0$ by \cite{mardesic} II.6.2
Lemma 2, therefore $\phi$ is a monomorphism. In particular there
exists $N\in\N$ big enough such that $\phi_{n,m_n}$ is an
injective homomorphism for every $n\geq N$. We set $N=1$ if
$\nu_\um{M}=\infty$. Now it is easy to see that the injection
$\psi_n\colon k\grupo{C_{n+1}}\hookrightarrow k\grupo{C_n}$ given
by the restriction of $\psi$ is the kernel of the next composition
whenever $n\geq N$
\begin{equation*}\label{unadetantas}
k\grupo{C_n}\st{\tau_0}\rightarrow
k\grupo{A_n}\twoheadrightarrow\frac{k\grupo{A_n}}
{\varphi(k\grupo{B_{m_n}})}=U^\varphi_{n,m_n}.
\end{equation*} Any morphism
$\upsilon\colon k\grupo{D}_\chi\rightarrow\um{C}_{\nu_\um{M}}$ is
represented by a controlled homomorphism $\upsilon_0\colon
k\grupo{D}_\chi\rightarrow k\grupo{C}_\gamma$. Suppose that
$\tau{\upsilon}=0$. This means that there exists another
controlled homomorphism $\eta\colon k\grupo{D}_\chi\rightarrow
k\grupo{B}_\beta$ with $\tau_0\upsilon_0=\varphi\eta$. By the
alternative characterization of controlled homomorphisms given in
Subsection \ref{finite} we see that there exists an increasing
sequence of natural numbers $\set{p_n}_{n\geq 1}$ such that
$\upsilon_0(k\grupo{D_{p_n}})\subset k\grupo{C_n}$ and
$\eta(k\grupo{D_{p_n}})\subset k\grupo{B_{m_n}}$, hence if $n\geq
N$ then there exists a unique homomorphism $\sigma_n\colon
k\grupo{D_{p_n}}\rightarrow k\grupo{C_{n+1}}$ such that
$\psi_n\sigma_n\colon k\grupo{D_{p_n}}\rightarrow k\grupo{C_n}$ is
the restriction of $\upsilon_0$. If $\sigma'\colon
k\grupo{{}_{p_{N}-1}D}\rightarrow k\grupo{C}$ is any homomorphism
such that $\psi\sigma'$ coincides with the restriction of
$\upsilon_0$ to $k\grupo{{}_{p_{N}-1}D}$ we define the controlled
homomorphism $\sigma\colon k\grupo{D}_\chi\rightarrow
k\grupo{C}_\gamma$ by $\sigma(d)=\sigma_n(d)$ if $d\in
{}_{p_{n+1}-1}D_{p_n}$ $(n\geq N)$ and $\sigma(d)=\sigma'(d)$ if
$d\in{}_{p_{N}-1}D$. This controlled homomorphism satisfies
$\psi\sigma=\upsilon_0$ therefore $\upsilon=0$ and hence $\tau$ is
a monomorphism.

Since $\tau$ is a monomorphism if we define $\um{N}=\coker\tau$ we
get by (\ref{medioin2}) that $\um{M}\simeq
\um{C}_{\nu_\um{M}}\oplus\um{N}$. Now one can check that
$\lambda_\um{N}=\mu_\um{N}=\nu_\um{N}=0$ by using that
$\um{N}=\coker[(\varphi,\tau_0)\colon k\grupo{B}_\beta\oplus
k\grupo{C}_\gamma\rightarrow k\grupo{A}_\alpha]$.
\end{proof}

\begin{lem}\label{es0}
If $\um{M}$ is a f. p. $\um{R}$-module with
$\lambda_\um{M}=\mu_\um{M}=\nu_\um{M}=0$ then $\um{M}=0$.
\end{lem}

\begin{proof}
If $\um{M}=\coker[\varphi\colon k\grupo{B}_\beta\rightarrow
k\grupo{A}_\alpha]$ the conditions of the statement are equivalent
to $\varphi(k\grupo{B})=k\grupo{A}$ and $\lim
U^\varphi_\bullet=0$. In the proof of (\ref{yvan2}) we checked
that $U^\varphi_\bullet$ satisfies the Mittag-Leffler property,
hence $U^\varphi_\bullet=0$ is a trivial pro-vector space by
\cite{mardesic} II.6.2 Lemma 2. This means that if
$\set{m_n}_{n\geq 1}$ is an increasing sequence such that
$\varphi(k\grupo{B_{m_n}})\subset k\grupo{A_n}$ (see Subsection
\ref{finite}) then there exists another increasing sequence
$\set{p_n}_{n\geq 1}$ such that the next bonding homomorphisms are
trivial
\begin{equation*}
U^\varphi_{p_{n+1},m_{p_{n+1}}}=\frac{k\grupo{A_{p_{n+1}}}}{\varphi(k\grupo{B_{m_{p_{n+1}}}})}
\st{0}\To\frac{k\grupo{A_{p_n}}}{\varphi(k\grupo{B_{m_{p_n}}})}
=U^\varphi_{p_n,m_{p_n}}.
\end{equation*} It is,
$k\grupo{A_{p_{n+1}}}\subset\varphi(k\grupo{B_{m_{p_n}}})$. Hence
we can define a controlled homomorphism $\psi\colon
k\grupo{A}_\alpha\rightarrow k\grupo{B}_\beta$ sending
$a\in{}_{p_{n+2}-1}A_{p_{n+1}}$ $(n\geq 1)$ to any element $b\in
B_{m_{p_n}}$ such that $\varphi(b)=a$, and if $a\in {}_{p_2}A$ we
take any $\psi(a)=b\in k\grupo{B}$ such that $\varphi(b)=a$. This
morphism satisfies $\varphi\psi=1$ hence $\varphi$ is an
epimorphism and $\um{M}=\coker\varphi=0$.
\end{proof}

\section{Representations of the $n$-subspace quiver}\label{n}

The \emph{$n$-subspace quiver} $Q_n$ is the next directed graph
$$\xygraph{*+=[o]+[F]{0} (:@{<-}[ull] *+=[o]+[F]{1}, :@{<-}[ul]
*+=[o]+[F]{2}, :@{}[u] {} :@{.}[r] {}, :@{<-}[urr]
*+=[o]+[F]{n})}$$ Fixed any field $k$, a \emph{representation}
$\ul{V}$ of $Q_n$ is a diagram of $k$-vector spaces indexed by
$Q_n$, it is, $n+1$ vector spaces $V_0, V_1, \dots V_n$ together
with homomorphisms $V_i\rightarrow V_0$  $(1\leq i\leq n)$.
Morphisms of representations are commutative diagrams. The
category $\rep_{Q_n}$ of representations of $Q_n$ is an abelian
category. It is equivalent to the category of $kQ_n$-modules,
where $kQ_n$ is the path algebra of $Q_n$, whose dimension is
$\dim kQ_n=2n+1$. A representation is said to be
finite-dimensional provided $V_i$ is a finite-dimensional vector
space for every $0\leq i\leq n$. Finitely presented (or
equivalently finite-dimensional) $kQ_n$-modules correspond to
finite-dimensional representations under this equivalence and
indecomposable representations are finite-dimensional, hence by
the classical Krull-Schmidt theorem the monoid
$\operatorname{Iso}(\rep^{\mathrm{fin}}_{Q_n})$ of isomorphisms
classes of finite-dimensional representations of $Q_n$ is the free
abelian monoid generated by the (isomorphism classes of)
indecomposable representations. The representation type of the
quiver $Q_n$ is that of its path algebra.

An \emph{$n$-subspace} $\ul{V}$ is a representation of $Q_n$ such
that the homomorphisms $V_i\rightarrow V_0$ are inclusions of
subspaces $V_i\subset V_0$ $(1\leq i\leq n)$. The category
$\sub_n$ (resp. $\sub^\mathrm{fin}_n$) of (finite-dimensional)
$n$-subspaces is a full additive subcategory of $\rep_{Q_n}$
(resp. $\rep^\mathrm{fin}_{Q_n}$). In fact direct summands in
$\rep_{Q_n}$ of $n$-subspaces are also $n$-subspaces, hence

\begin{prop}
$\operatorname{Iso}(\sub^\mathrm{fin}_n)$ is the free abelian
monoid generated by the isomorphism classes of indecomposable
$n$-subspaces.
\end{prop}

Up to isomorphism there are just $n$ indecomposable
representations of $Q_n$ which are not $n$-subspaces, namely those
with $V_i=k$ for some $1\leq i\leq n$ and $V_j=0$ if $j\neq i$,
therefore

\begin{prop}
The representation type of $\sub^\mathrm{fin}_n$ is the same as
the $n$-subspace quiver.
\end{prop}

We say that an $n$-subspace $\ul{V}$ is \emph{rigid} provided
$V_i\subset \sum_{j\neq 0,i}V_j$ $(1\leq i\leq n)$ and
$V_0=\sum_{i=1}^nV_i$. As before, the category
$\sub_n^\mathrm{rig}$ (resp. $\sub_n^\mathrm{fr}$) of
(finite-dimensional) rigid $n$-subspaces is an additive (small)
subcategory of $\sub_n$ (resp. $\sub^\mathrm{fin}_n$) and direct
summands of rigid $n$-subspaces are also rigid, so

\begin{prop}
$\operatorname{Iso}(\sub^\mathrm{fr}_n)$ is the free abelian
monoid generated by the isomorphism classes of indecomposable
rigid $n$-subspaces.
\end{prop}

There is an additive ``rigidification'' functor
\begin{equation*}
\sub_n\To\sub_n^\mathrm{rig}\colon \ul{V}\mapsto
\ul{V}^\mathrm{rig}
\end{equation*} given by
$V^\mathrm{rig}_i=V_i\cap\left(\sum_{j\neq 0,i}V_j\right)$ $(1\leq
i\leq n)$ and $V_0^\mathrm{rig}=\sum_{i=1}^nV^\mathrm{rig}_i$,
which is right-adjoint to the inclusion
$\sub_n^\mathrm{rig}\subset\sub_n$ and preserves
finite-dimensional objects. The unit of this adjunction is the
obvious natural inclusion $\ul{V}^\mathrm{rig}\subset\ul{V}$,
which is an equality if and only if $\ul{V}$ is already rigid.

In order to determine indecomposable rigid $n$-subspaces we
consider the full inclusions of additive categories
$\uf{F}^i\colon\sub_1\rightarrow\sub_n$ $(1\leq i\leq n)$ sending
a 1-subspace $\ul{W}$ to the $n$-subspace
$\uf{F}^i\ul{W}={}^i\ul{W}$ with ${}^iW_0=W_0$, ${}^i{W}_i=W_1$
and ${}^i{W}_j=0$ otherwise.

\begin{prop}\label{rigiP}
The natural inclusion $\ul{V}^\mathrm{rig}\subset \ul{V}$ admits a
(not natural) retraction in $\sub_n$. More precisely, there exist
1-subspaces $\ul{V}^i$ $(1\leq i\leq n)$ and an isomorphism
$\ul{V}\simeq (\bigoplus_{i=1}^n\uf{F}^i\ul{V}^i)\oplus
\ul{V}^\mathrm{rig}$ such that the natural inclusion
$\ul{V}^\mathrm{rig}\subset \ul{V}$ corresponds to the inclusion
of the direct summand.
\end{prop}

\begin{proof}
By using the definition of $\ul{V}^\mathrm{rig}$ we see that there
is a short exact sequence of vector spaces
$$\bigoplus_{i=1}^n\frac{V_i}{V^\mathrm{rig}_i}\hookrightarrow
\frac{V_0}{V_0^\mathrm{rig}}\twoheadrightarrow\frac{V_0}{\sum_{i=1}^nV_i}.$$
We define the $1$-subspaces $\ul{V}^i$ $(1\leq i\leq n)$ as
$V^1_0=\frac{V_0}{\sum_{i=1}^nV_i}\oplus\frac{V_1}{V^\mathrm{rig}_1}$,
$V^1_1=\frac{V_1}{V^\mathrm{rig}_1}$ and
$V^i_0=V^i_1=\frac{V_i}{V^\mathrm{rig}_i}$ if $1<i\leq n$. Now the
isomorphism of the statement follows from the former exact
sequence.
\end{proof}

By this proposition an indecomposable $n$-subspace is rigid unless
it is isomorphic to $\uf{F}^i\ul{V}$ for some indecomposable
1-subspace $\ul{V}$ and $1\leq i\leq n$. It is known that such
$\ul{V}$ must be either $k\rightarrow k$ or $0\rightarrow k$, so
there are just $2n$ indecomposable $n$-subspaces which are not
rigid, and $3n$ indecomposable representations of $Q_n$ which are
not rigid $n$-subspaces, in particular

\begin{prop}\label{esigual}
The category $\sub^\mathrm{fr}_n$ has the same representation type
as the $n$-subspace quiver.
\end{prop}

\begin{rem}\label{sesabe}
The representation type of the $n$-subspace quiver is well-known.
It is finite for $n<4$, tame for $n=4$ and wild if $n>4$, see
\cite{gabriel} and \cite{nazarova}.
\end{rem}

In \cite{gabriel} the finite sets of indecomposable
representations of $Q_n$ are described for $n<4$ hence discarding
the $3n$ indecomposable representations previously described which
are not rigid $n$-subspaces we get the next result.

\begin{prop}\label{losn}
The following are complete lists of (representatives of the
isomorphism classes of) indecomposable rigid $n$-subspaces for
$n<4$
\begin{itemize}
\item $n=1$, none, \item $n=2$, $\ul{V}^{(2,1)}=\left(k\rightarrow
k\leftarrow k\right)$, \item $n=3$,
$$\ul{V}^{(3,1)}=\left(\begin{array}{ccccc}
   &  & k &  &  \\
   &  & \downarrow &  &  \\
  k & \rightarrow & k & \leftarrow & 0
\end{array}\right),\;\;\;
\ul{V}^{(3,2)}=\left(\begin{array}{ccccc}
   &  & 0 &  &  \\
   &  & \downarrow &  &  \\
  k & \rightarrow & k & \leftarrow & k
\end{array}\right),$$
$$\ul{V}^{(3,3)}=\left(\begin{array}{ccccc}
   &  & k &  &  \\
   &  & \downarrow &  &  \\
  0 & \rightarrow & k & \leftarrow & k
\end{array}\right),\;\;\;
\ul{V}^{(3,4)}=\left(\begin{array}{ccccc}
   &  & k &  &  \\
   &  & \downarrow &  &  \\
  k & \rightarrow & k & \leftarrow & k
\end{array}\right),$$
$$\ul{V}^{(3,5)}=\left(\begin{array}{ccccc}
   &  & k\grupo{x+y} &  &  \\
   &  & \downarrow &  &  \\
  k\grupo{x} & \rightarrow & k\grupo{x,y} & \leftarrow & k\grupo{y}
\end{array}\right).$$
\end{itemize}
\end{prop}


\begin{rem}\label{los4}
In \cite{nazarova} there is a (not finite) list of
indecomposable representations of $Q_4$ 
. We do not include the list here
because it is quite tedious to describe, however the interested
reader can easily find and remove the 12 indecomposable
representations of $Q_4$ which are not rigid $4$-subspaces,
obtaining in this way a complete list of indecomposable rigid
$4$-subspaces 
.
\end{rem}

\section{Finitely presented $k(\n{n})$-modules and finite-dimensional $n$-subspaces}\label{modandn}

Given an $n$-subspace $\ul{V}$ we define $\uf{M}\ul{V}\colon
\M_k(\ol{T}_n)^{op}\rightarrow\Ab$ as the additive functor which
sends an object $k\grupo{A}_\alpha$ to the vector subspace
$\uf{M}\ul{V}(k\grupo{A}_\alpha)\subset \hom_k(k\grupo{A},V_0)$
formed by the homomorphisms $\phi\colon k\grupo{A}\rightarrow V_0$
such that there exists $M\geq 1$ depending on $\phi$ satisfying
$\phi(A^i_M)\subset V_i$ $(1\leq i\leq n)$. This construction
defines an exact full inclusion of additive categories
\begin{equation}\label{amod}
\uf{M}\colon\sub_n\rightarrow\mo(\M_k(\ol{T}_n)).
\end{equation}

\begin{prop}\label{esfp}
If $\ul{V}$ is a finite-dimensional $n$-subspace then the
$\M_k(\ol{T}_n)$-module $\uf{M}\ul{V}$ is finitely presented.
\end{prop}

\begin{proof}
Let $\set{w_1,\dots w_d}$ be a basis of $V_0$, $\set{w^i_1,\dots
w^i_{d_i}}$ a basis of $V_i$ $(1\leq i\leq n)$, and $\phi_i\colon
V_i\rightarrow V_0$ the inclusion. We define the sets
$D=\set{{}_mw^i_1,\dots {}_mw^i_{d_i}\,;\,1\leq i\leq n,m\geq 1}$
and $C=D\sqcup \set{w_1,\dots w_d}$, and the height functions
$\gamma\colon C\rightarrow T^0_n$  and $\delta\colon D\rightarrow
T^0_n$ with $\gamma(w_j)=v_0$ $(1\leq j\leq d)$  and
$\gamma({}_mw^i_j)=\delta({}_mw^i_j)=v^i_m$ $(1\leq i\leq n, 1\leq
j\leq d_i,m\geq 1 )$. Let $\rho\colon k\grupo{D}_\delta\rightarrow
k\grupo{C}_\gamma$ be the controlled homomorphism defined as
$\rho({}_mw^i_j)={}_mw^i_j-{}_{m-1}w^i_j$ if $m>1$ and
$\rho({}_1w^i_j)={}_1w^i_j-\phi_i(w^i_j)$ otherwise, and $p\colon
k\grupo{C}_\gamma\rightarrow\uf{M}\ul{V}$ the
$\M_k(\ol{T}_n)$-module morphism determined by the $k$-vector
space homomorphism $p_0\colon k\grupo{C}\rightarrow V_0$ with
$p_0(w_j)=w_j$ and $p_0({}_mw^i_j)=\phi_i(w^i_j)$. Here we use
that
$\hom(k\grupo{C}_\gamma,\uf{M}\ul{V})=\uf{M}\ul{V}(k\grupo{C}_\gamma)$
by Yoneda's lemma. Now it is immediate to check that
$\uf{M}\ul{V}=\coker\rho$ and $p$ is the natural projection.
\end{proof}

One can check by using the finite presentation constructed in the
proof of the former proposition that

\begin{cor}\label{seanula}
If $\ul{V}$ is a finite-dimensional rigid $n$-subspace then
$\Phi_n([\uf{M}\ul{V}])=0$.
\end{cor}

By (\ref{esfp}) the additive functor in (\ref{amod}) restricts to
a functor
$\uf{M}\colon\sub_n^\mathrm{fin}\rightarrow\fp(\M_k(\ol{T}_n))$.
Now we are going to construct a functor in the opposite direction.
For this if $\varphi\colon k\grupo{B}_\beta\rightarrow
k\grupo{A}_\alpha$ is a morphism in $\M_k(\ol{T}_n)$ we define the
$k$-vector spaces $(1\leq i\leq n)$
$$W^\varphi_i=\frac{\bigcap_{m\geq 1}\left\{\left[k\grupo{A^i_m}+
\varphi(k\grupo{B})\right]\cap\left[\sum_{j\neq i} k\grupo{A_m^j}+
\varphi(k\grupo{B})\right]\right\}}{\varphi(k\grupo{B})}.$$

\begin{prop}\label{finitud}
The vector space $W_i^\varphi$ is finite-dimensional $(1\leq i\leq
n)$.
\end{prop}

This proposition is an immediate consequence of the next

\begin{lem}\label{finilem}
For any $m\geq 1$ $$\dim \frac{\left[ k\grupo{A^i_m}+
\varphi(k\grupo{B})\right]\cap\left[\sum_{j\neq i} k\grupo{A_m^j}+
\varphi(k\grupo{B})\right]}{\varphi(k\grupo{B})}<\aleph_0.$$
\end{lem}

\begin{proof}\renewcommand{\theequation}{\alph{equation}}\setcounter{equation}{0}
One readily checks that $$\left[ k\grupo{A^i_m}+
\varphi(k\grupo{B})\right]\cap\left[\sum_{j\neq i} k\grupo{A_m^j}+
\varphi(k\grupo{B})\right]=$$$$k\grupo{A^i_m}\cap\left[\sum_{j\neq
i} k\grupo{A_m^j}+
\varphi(k\grupo{B})\right]+\varphi(k\grupo{B}),$$ therefore
\begin{equation}\label{una}
\frac{\left[ k\grupo{A^i_m}+
\varphi(k\grupo{B})\right]\cap\left[\sum_{j\neq i} k\grupo{A_m^j}+
\varphi(k\grupo{B})\right]}{\varphi(k\grupo{B})}\simeq
\end{equation}
$$\frac{k\grupo{A^i_m}\cap\left[\sum_{j\neq i} k\grupo{A_m^j}+
\varphi(k\grupo{B})\right]}{k\grupo{A^i_m}\cap\varphi(k\grupo{B})}.$$

Since $\varphi$ is a controlled homomorphism there exists $M\geq
1$ such that $\varphi(B_M^i)\subset k\grupo{A^i_m}$ $(1\leq i\leq
n)$, hence
\begin{equation}\label{otra}
\sum_{j\neq i}k\grupo{A^j_m}+\varphi(k\grupo{B})= \sum_{j\neq
i}k\grupo{A^j_m}+\varphi(k\grupo{{}_{M-1}B})+\varphi(k\grupo{B^i_M}).
\end{equation}

The set ${}_{M-1}B$ is finite, hence there exists $N\geq 0$ big
enough with $\varphi({}_{M-1}B)\subset k\grupo{{}_NA}$. Let us
check that the next homomorphism induced by the inclusion
${}_NA^i_m\subset A^i_m$ is an isomorphism
\begin{equation}\label{reotra}
\frac{k\grupo{{}_NA^i_m}\cap\left[\sum_{j\neq i} k\grupo{A_m^j}+
\varphi(k\grupo{B})\right]}{k\grupo{{}_NA^i_m}\cap\varphi(k\grupo{B})}\rightarrow\frac{k\grupo{A^i_m}\cap\left[\sum_{j\neq
i} k\grupo{A_m^j}+
\varphi(k\grupo{B})\right]}{k\grupo{A^i_m}\cap\varphi(k\grupo{B})}.
\end{equation}
The injectivity is obvious. Now by (\ref{otra}) an arbitrary
element in the range of (\ref{reotra}) is represented by an
element $a_i\in k\grupo{A^i_m}$ such that there are $a_j\in
k\grupo{A^j_m}$ $(j\neq i)$, $a'_i\in k\grupo{A^i_m}\cap
\varphi(k\grupo{B})$, $b\in \varphi(k\grupo{{}_{M-1}B})$, and
$b_i\in \varphi(k\grupo{B^i_M})$ such that $$\sum_{j\neq
i}a_j+b+b_i=a_i+a_i',$$ but $a_i+a_i'-b_i\in k\grupo{A^i_m}$,
$(\oplus_{j=1}^nk\grupo{A^j_m})\cap
k\grupo{{}_NA}=\oplus_{j=1}^nk\grupo{{}_NA^j_m}$, and $\sum_{j\neq
i}a_j-(a_i+a_i'-b_i)=b\in k\grupo{{}_NA}$, therefore
$a_j,a_i+a_i'-b_i\in k\grupo{{}_NA}$ $(j\neq i)$, and $a_i'-b_i\in
k\grupo{A^i_m}\cap\varphi\grupo{B}$, so $a_i+a_i'-b_i$ represents
the same element as $a_i$ in the range of (\ref{reotra}), hence
the homomorphism (\ref{reotra}) is surjective. Now the proposition
follows from the isomorphisms (\ref{una}) and (\ref{reotra}), and
the finiteness of the set ${}_NA^i_m$.
\end{proof}

By (\ref{finitud}) the vector space $W_0^\varphi=\sum_{i=1}^n
W_i^\varphi$ together with the subspaces $W^\varphi_1, \dots
W^\varphi_n$ define a finite dimensional $n$-subspace
$\ul{W}^\varphi$.

\begin{prop}\label{s}
There is an additive functor
$\uf{S}\colon\fp({\M_k(\ol{T}_n)})\rightarrow\sub_n^\mathrm{fin}$
which sends $\um{M}=\coker\varphi$ to
$\uf{S}\um{M}=\ul{W}^\varphi$.
\end{prop}

\begin{proof}
By using the alternative description of controlled homomorphisms
in $\M_k(\ol{T}_n)$ given in Subsection \ref{finite} one can
easily check hat the correspondence $\varphi\mapsto
\ul{W}^\varphi$ is a functor from $\pair({\M_k(\ol{T}_n)})$ to the
category of $n$-subspaces. Furthermore, this functor factors
through the natural equivalence relation $\sim$, therefore the
proposition follows by (\ref{cokereq}).
\end{proof}

The functors $\uf{S}$ and $\uf{M}$ are not adjoint. Moreover, one
readily checks that

\begin{prop}
The $n$-subspace $\uf{S}\um{M}$ is rigid for every f. p.
$\M_k(\ol{T}_n)$-module $\um{M}$. Moreover, given a
finite-dimensional $n$-subspace $\ul{V}$ there is a natural
isomorphism $\uf{S}\uf{M}\ul{V}\simeq{\ul{V}^\mathrm{rig}}$.
\end{prop}

For the second statement of the former proposition one uses the
finite presentations constructed in the proof of (\ref{esfp}).

\begin{cor}\label{qei}
The image of the functor $\uf{S}$ is the category of
finite-dimensional rigid $n$-subspaces.
\end{cor}

\begin{cor}\label{qei2}
For every f. p. $\M_k(\ol{T}_n)$-module $\um{M}$ there is a
natural isomorphism $\uf{SMS}\um{M}\simeq\uf{S}\um{M}$.
\end{cor}

As we pointed out in the introduction a key step to obtain a
presentation of $\operatorname{Iso}(\fp(\M_k(\ol{T}_n)))$ is
relating the decomposition problem in $\fp(\M_k(\ol{T}_n))$ to the
decomposition problem in $\fp(\M_k(\ol{T}_1))$ and
$\sub^{\mathrm{fin}}_n$. This is what we do in the next two
propositions.

\begin{prop}\label{tampoco}
For any f. p. $\M_k(\ol{T}_n)$-module $\um{M}$ there exists
another one $\um{N}$ with $\uf{S}\um{N}=0$ such that
$\um{M}\simeq\um{N}\oplus\uf{M}\uf{S}\um{M}$.
\end{prop}

\begin{proof}\renewcommand{\theequation}{\alph{equation}}\setcounter{equation}{0}
Suppose that $\um{M}$ is the cokernel of $\varphi\colon
k\grupo{B}_\beta\rightarrow k\grupo{A}_\alpha$ in
$\M_k(\ol{T}_n)$. By the alternative description of controlled
homomorphisms given in Subsection \ref{finite} we can choose an
increasing sequence of natural numbers $\set{M_m}_{m\geq 1}$ with
$\varphi(B^i_{M_m})\subset k\grupo{A^i_m}$ $(1\leq i\leq n)$. We
define the inverse systems of vector spaces $X^i_\bullet,
Y^i_\bullet, Z^i_\bullet$ $(1\leq i\leq n)$ indexed by $\N$ in the
following way
$$X^i_m=\frac{\varphi(k\grupo{B})}{\varphi(k\grupo{B^i_{M_m}})},$$
$$Y^i_m=\frac{\left[ k\grupo{A^i_m}+
\varphi(k\grupo{B})\right]\cap\left[\sum_{j\neq i} k\grupo{A_m^j}+
\varphi(k\grupo{B})\right]}{\varphi(k\grupo{B^i_{M_m}})},$$
$$Z^i_m=\frac{\left[ k\grupo{A^i_m}+
\varphi(k\grupo{B})\right]\cap\left[\sum_{j\neq i} k\grupo{A_m^j}+
\varphi(k\grupo{B})\right]}{\varphi(k\grupo{B})}.$$ The bonding
homomorphisms are induced by the obvious inclusions of vector
spaces. The short exact sequences $X_m^i\hookrightarrow
Y^i_m\st{p_i^m}\twoheadrightarrow Z^i_m$ are compatible with the
bonding homomorphisms, so they give rise to short exact sequences
$X_\bullet^i\hookrightarrow Y^i_\bullet\st{p_i}\twoheadrightarrow
Z^i_\bullet$ in the abelian pro-category of pro-vector spaces.
Moreover, $\lim Z^i_\bullet=\cap_{m\geq 1}Z^i_m=W^\varphi_i$. Let
$\psi_i\colon W^\varphi_i\rightarrow Z^i_\bullet$ be the canonical
pro-morphism, which is induced by the inclusions
$W^\varphi_i\subset Z^i_m$. Here we regard $W^\varphi_i$ as an
inverse system indexed by a singleton.

The bonding homomorphisms of the inverse system $X^i_\bullet$ are
surjective, therefore ${\lim }^1X^i_\bullet=0$, and by (\ref{gro})
$\ext^1(W^\varphi_i,X^i_\bullet)=0$, so there exists a
pro-morphism $\tilde{\psi}_i$ such that the next diagram commutes
$$\xymatrix{&&W^\varphi_i\ar[d]^{\psi_i}\ar[ld]_{\tilde{\psi}_i}\\X_\bullet^i\ar@{^{(}->}[r]&
Y^i_\bullet\ar@{->>}[r]^{p_i}& Z^i_\bullet}$$ The pro-morphism
$\tilde{\psi}_i$ is represented by a sequence of homomorphisms
$\tilde{\psi}_i^m\colon W^\varphi_i\rightarrow Y^i_m$ $(m\geq 1)$
which are compatible with the bonding homomorphisms of
$Y^i_\bullet$, and such that the composition
$p^m_i\tilde{\psi}^m_i\colon W^\varphi_i\subset Z^i_m$ is the
inclusion.

If $\set{a_1^i,\dots a^i_{d_i}}$ is a basis of  $W^\varphi_i$ we
can choose elements $\set{{}_ma_1^i,\dots {}_ma_{d_i}^i}\subset
k\grupo{A^i_m}$ $(m\geq 1)$ such that
$\tilde{\psi}_i^m(a_j^i)={}_ma^i_j+\varphi(k\grupo{B^i_{M_m}})$.
In particular, since the homomorphisms $\tilde{\psi}_i^m$ are
compatible with the bonding homomorphisms of $Y^i_\bullet$, we see
that there are elements ${}_{m+1}b^i_j\in k\grupo{B^i_{M_m}}$
satisfying ${}_{m+1}a^i_j-{}_ma^i_j=\varphi({}_{m+1}b^i_j)$.
Moreover, let $\set{a_1,\dots a_d}$ be a basis of $W^\varphi_0$, 
$\sigma\colon k\grupo{A}/\varphi(k\grupo{B})\hookrightarrow
k\grupo{A}$ a splitting of the natural projection, and elements
${}_1b_j^i\in k\grupo{B}$ such that
$\varphi({}_1b_j^i)={}_1a_j^i-\sigma({}_1a_j^i+\varphi(k\grupo{B}))$.

If $\rho$ is the finite presentation of $\uf{M}\uf{S}\um{M}$
constructed in the proof of (\ref{esfp}), there is a morphism
$\tau\colon\rho\rightarrow\varphi$ in $\pair({\M_k(\ol{T}_n)})$
given by $\tau_0(w_i)=\sigma(a_i)$, $\tau_0({}_mw^i_j)={}_ma^i_j$,
and $\tau_1({}_mw^i_j)={}_mb^i_j$. This morphism induces a
$\M_k(\ol{T}_n)$-module morphism
$\tau\colon\uf{M}\uf{S}\um{M}\rightarrow\um{M}$. Now we are going
to construct a retraction of $\tau$.

By (\ref{finilem}) $Z^i_m$ is always finite-dimensional and
$W^\varphi_i=\cap_{m\geq 1}Z^i_m$, hence there exists $N\geq 1$
such that $W^\varphi_i=Z^i_N$ for every $1\leq i\leq n$. Let
$\ul{V}$ be the $n$-subspace given by
$V_0=k\grupo{A}/\varphi(k\grupo{B})$ and
$V_i=\left[k\grupo{A^i_N}+\varphi(k\grupo{B})\right]/\varphi(k\grupo{B})$.
Clearly $\uf{S}\um{M}=\ul{V}^\mathrm{rig}\subset \ul{V}$, hence by
(\ref{rigiP}) there is a retraction $r\colon
\ul{V}\rightarrow\uf{S}\um{M}$. By Yoneda's lemma
$\hom_\mathcal{R}(k\grupo{A}_\alpha,\uf{M}\ul{V})=\uf{M}\ul{V}(k\grupo{A}_\alpha)$.
The natural projection $k\grupo{A}\twoheadrightarrow V_0$ give
rise to a $\M_k(\ol{T}_n)$-module morphism $\upsilon_0\colon
k\grupo{A}_\alpha\rightarrow\uf{M}\ul{V}$ such that
$\upsilon_0\varphi=0$. Moreover, since $\um{M}=\coker\varphi$ then
$\hom_\mathcal{R}(\um{M},\uf{M}\ul{V})=\ker\hom_\mathcal{R}(\varphi,\uf{M}\ul{V})$,
in particular $\upsilon_0$ determines a $\M_k(\ol{T}_n)$-module
morphism $\upsilon\colon\um{M}\rightarrow\uf{M}\ul{V}$. One
readily checks that the composition $(\uf{M}r)\upsilon_0\tau_0$
coincides with the natural projection $p\colon
k\grupo{C}_\gamma\twoheadrightarrow\uf{M}\uf{S}\um{M}=\coker\rho$
defined in the proof of (\ref{esfp}), hence
$(\uf{M}r)\upsilon\tau=1$ is the identity on $\uf{M}\uf{S}\um{M}$,
and $(\uf{M}r)\upsilon$ is the desired retraction of $\tau$. Now
if we take $\um{N}$ to be the cokernel of $\tau$ the proposition
follows since $\um{M}=\um{N}\oplus\uf{M}\uf{S}\um{M}$, by
(\ref{qei2})
$\uf{S}\um{M}=\uf{S}\um{N}\oplus\uf{S}\uf{M}\uf{S}\um{M}\simeq
\uf{S}\um{N}\oplus\uf{S}\um{M}$ and hence $\uf{S}\um{N}=0$. Here
we use that the monoid $\operatorname{Iso}(\sub^\mathrm{fin}_n)$
is free and hence cancelative, compare Section \ref{n}.
\end{proof}


In the next proposition we use the change of coefficients
$\uf{F}^i_*$ associated to the additive functors
$\uf{F}^i\colon\M_k(\ol{T}_1)\rightarrow\M_k(\ol{T}_n)$ in Remark
\ref{teincluyo}.

\begin{prop}\label{nuse}
Given a f. p. $\M_k(\ol{T}_n)$-module $\um{M}$, $\uf{S}\um{M}=0$
if and only if there exist f. p. $\M_k(\ol{T}_1)$-modules
$\um{M}_i$ $(1\leq i\leq n)$ with $\um{M}\simeq
\uf{F}^1_*\um{M}_1\oplus\cdots\oplus\uf{F}^n_*\um{M}_n$.
\end{prop}

\begin{proof}\renewcommand{\theequation}{\alph{equation}}\setcounter{equation}{0}
It is easy to see that $\uf{S}\uf{F}^i_*=0$ $(1\leq i\leq n )$,
and $\uf{S}$ is additive, so the implication $\Leftarrow$ follows.
Now suppose that $\um{M}=\coker\left[\varphi\colon
k\grupo{B}_\beta\rightarrow k\grupo{A}_\alpha\right]$ and
$\uf{S}\um{M}=0$. Since finite-dimensional vector spaces are
artinian, by (\ref{finilem}) there exists $m\geq 1$ big enough
such that for every $1\leq i\leq n$,
$$\frac{\left[ k\grupo{A^i_m}+
\varphi(k\grupo{B})\right]\cap\left[\sum_{j\neq i} k\grupo{A_m^j}+
\varphi(k\grupo{B})\right]}{\varphi(k\grupo{B})}=0,$$ it is, the
following equality holds (the isomorphism on the right always
holds)
$$\frac{\sum_{i=1}^n k\grupo{A_m^i}+
\varphi(k\grupo{B})}{\varphi(k\grupo{B})}=\bigoplus_{i=1}^n\frac{k\grupo{A_m^i}+
\varphi(k\grupo{B})}{\varphi(k\grupo{B})}\simeq
\bigoplus_{i=1}^n\frac{k\grupo{A_m^i}}{k\grupo{A_m^i}\cap\varphi(k\grupo{B})}.$$
This is equivalent to state that
\begin{equation}\label{princi}
\left[\bigoplus_{i=1}^nk\grupo{A_m^i}\right]\cap\varphi(k\grupo{B})=\bigoplus_{i=1}^n
\left[k\grupo{A_m^i}\cap\varphi(k\grupo{B})\right].
\end{equation}

By the characterization of controlled homomorphisms in Subsection
\ref{finite} there exists $M\geq 1$ with $\varphi(B^i_M)\subset
k\grupo{A^i_m}$ $(1\leq i\leq n)$. Let $K$ be the kernel of the
vector space homomorphism underlying to $\varphi$, it is
$K=\varphi^{-1}(0)$. There is a finite set $\set{b_1,\dots
b_d}\subset k\grupo{B}$ which projects to a basis of
$$\frac{K+\left(\bigoplus_{i=0}^nk\grupo{B^i_M}\right)}
{\bigoplus_{i=0}^nk\grupo{B^i_M}}\simeq \frac{K}
{K\cap\left(\bigoplus_{i=0}^nk\grupo{B^i_M}\right)},$$ since this
vector space is contained in
\begin{equation*}
\frac{k\grupo{B}}{\bigoplus_{i=0}^nk\grupo{B^i_M}}\simeq
k\grupo{{}_{M-1}B},
\end{equation*}
and ${}_{M-1}B$ is finite.

There is also a finite set $\set{a^i_1,\dots a^i_{d_i}}\subset
k\grupo{B}$ which projects to a basis of
\begin{equation*}
\frac{\varphi^{-1}\left(k\grupo{A^i_m}\cap
\varphi(k\grupo{B})\right)}{k\grupo{B^i_M}+K},
\end{equation*}
because $\varphi$ induces an isomorphism
\begin{equation*}
\frac{\varphi^{-1}\left(k\grupo{A^i_m}\cap
\varphi(k\grupo{B})\right)}{k\grupo{B^i_M}+K}\simeq
\frac{k\grupo{A^i_m}\cap
\varphi(k\grupo{B})}{\varphi(k\grupo{B^i_M})},
\end{equation*}
and always
\begin{equation}\label{inter}
k\grupo{A^i_m}\cap\left(\sum_{j\neq i}k\grupo{A^j_m}\right)=0,
\end{equation}
so
$k\grupo{A^i_m}\cap\left[\sum_{j=1}^n\varphi(k\grupo{B^j_M})\right]=\varphi(k\grupo{B^i_M})$,
and hence
\begin{equation*}
\frac{k\grupo{A^i_m}\cap
\varphi(k\grupo{B})}{\varphi(k\grupo{B^i_M})}=\frac{k\grupo{A^i_m}\cap
\varphi(k\grupo{B})}
{k\grupo{A^i_m}\cap\left[\sum_{j=1}^n\varphi(k\grupo{B^j_M})\right]}
\subset\frac{\varphi(k\grupo{B})}{\sum_{i=1}^n\varphi(k\grupo{B^i_M})}\simeq
\varphi(k\grupo{{}_{M-1}B}).
\end{equation*}

By (\ref{inter}) we have that
\begin{equation*}
\varphi^{-1}\left(k\grupo{A^i_m}\cap
\varphi(k\grupo{B})\right)\cap\left[\sum_{j\neq
i}\varphi^{-1}\left(k\grupo{A^j_m}\cap
\varphi(k\grupo{B})\right)\right]=K,
\end{equation*} therefore the set $\left[\sqcup_{i=1}^n\left(B^i_M
\sqcup\set{a^i_j}_{j=1}^{d_i}\right)\right]\sqcup\set{b_i}_{i=1}^d$
is linearly independent in $k\grupo{B}$. Moreover, it is a basis
of $\sum_{i=1}^n\varphi^{-1}\left(k\grupo{A^i_m}\cap
\varphi(k\grupo{B})\right)$, so in order to complete it to a basis
$\ul{B}$ of $k\grupo{B}$ we only need to add a finite set
$\set{b'_1,\dots b'_{d'}}\subset k\grupo{B}$ which projects to a
basis of the following vector space
$$\frac{k\grupo{B}}{\sum_{i=1}^n\varphi^{-1}\left(k\grupo{A^i_m}\cap
\varphi(k\grupo{B})\right)}.$$ This vector space is isomorphic to
\begin{equation}\label{basea}
\frac{\varphi(k\grupo{B})}{\bigoplus_{i=1}^nk\grupo{A^i_m}\cap
\varphi(k\grupo{B})}\subset\frac{k\grupo{A}}{\bigoplus_{i=1}^nk\grupo{A^i_m}}\simeq
k\grupo{{}_{m-1}A},
\end{equation}
and hence finite-dimensional. The inclusion (\ref{basea}) follows
from (\ref{princi}). Let $\set{a_1,\dots a_e}\subset k\grupo{A}$
be a basis of
\begin{equation*}
\frac{k\grupo{A}}{\varphi(k\grupo{B})+\left(\bigoplus_{i=1}^nk\grupo{A^i_m}\right)}.
\end{equation*}
By (\ref{basea}) $\ul{A}=\left(\sqcup_{i=1}^nA^i_m
\right)\sqcup\set{\varphi(b_i')}_{i=1}^{d'}\sqcup\set{a_i}_{i=1}^e$
is a basis of $k\grupo{A}$. Let $\ul{\alpha}\colon
\ul{A}\rightarrow T^0_n$, $\ul{\beta}\colon \ul{B}\rightarrow
T^0_n$ be the height functions defined as $\alpha$ and $\beta$
over $\sqcup_{i=1}^nA^i_m$ and $\sqcup_{i=1}^nB^i_M$ respectively,
and constant $v_0$ on the other elements. The identities
$k\grupo{A}=k\grupo{\ul{A}}$ and $k\grupo{B}=k\grupo{\ul{B}}$
induce controlled isomorphisms $\phi_1\colon
k\grupo{A}_\alpha\simeq k\grupo{\ul{A}}_{\ul{\alpha}}$ and
$\phi_2\colon k\grupo{B}_\beta\simeq
k\grupo{\ul{B}}_{\ul{\beta}}$, so if we define
$\psi=\phi_1\varphi\phi_2^{-1}$ then $\um{M}\simeq\coker\psi$. But
if we define the sets
${}^1\ul{A}=A^1_m\sqcup\set{\varphi(b_i')}_{i=1}^{d'}\sqcup\set{a_i}_{i=1}^e$,
${}^i\ul{A}=A^i_m$ $(1<i\leq n)$, ${}^1\ul{B}=B^1_M
\sqcup\set{a^1_j}_{j=1}^{d_1}\sqcup\set{b_i}_{i=1}^d\sqcup\set{b_i'}_{i=1}^{d'}$,
${}^i\ul{B}=B^i_M \sqcup\set{a^i_j}_{j=1}^{d_i}$ $(1<i\leq n)$,
and the height functions ${}^i\ul{\alpha}$ and ${}^i\ul{\beta}$ as
the restriction of $\ul{\alpha}$ and $\ul{\beta}$ to ${}^i\ul{A}$
and ${}^i\ul{B}$ respectively $(1\leq i\leq n)$, then we observe
that $(1\leq i\leq n)$
\begin{equation*}\label{tela1}
{}^i\ul{\alpha}({}^i\ul{A}), {}^i\ul{\beta}({}^i\ul{B}) \subset
\set{v_0}\cup\set{v^i_m}_{m\geq 1}\subset T^0_n,
\end{equation*}
\begin{equation*}\label{tela3}
\psi({}^i\ul{A})\subset k\grupo{{}^i\ul{B}},
\end{equation*}
and the proposition follows.
\end{proof}



\section{Classification of finitely presented
$k(\n{n})$-modules}\label{oju}

In this final section we complete the proofs of Theorems \ref{rep}
and \ref{generel}. For this, the crucial result is the next
theorem where we compute the monoid
$\operatorname{Iso}(\fp(k(\n{n})))$ in terms of the free abelian
monoid $\operatorname{Iso}(\sub^{\mathrm{fr}}_n)$.

\begin{thm}\label{jander}
The following monoid morphism is an isomorphism for every
$n\in\N$:
\begin{equation*}
\left(\Phi_n,\operatorname{Iso}(\uf{S})\right)\colon\operatorname{Iso}(\fp(k(\n{n})))
\st{\simeq}\To\N_{\infty,n}\times\prod_{i=1}^n\N_{\infty}
\times\prod_{i=1}^n\N_{\infty}\times\operatorname{Iso}(\sub^{\mathrm{fr}}_n).
\end{equation*}
\end{thm}

This theorem follows from the strongest results previously proven
in this paper. More precisely, the surjectivity of
$\left(\Phi_n,\operatorname{Iso}(\uf{S})\right)$ is a consequence
of (\ref{metiendo}), (\ref{imagenes}), (\ref{seanula}),
(\ref{qei}), (\ref{qei2}) and (\ref{nuse}). Furthermore, this
morphism is injective by the next theorem. In order to state it we
introduce the following notation. Given $d\in\N_{\infty,n}$ the
$k(\n{n})$-module $\um{A}_d$ is $\uf{F}^1_*\um{A}_d$ provided
$d\in\NN$ and $\oplus_{i\in S}\uf{F}^i_*\um{R}$ if $d=\infty_S$
for some $\emptyset\neq S\subset\{1,\dots n\}$.

\begin{thm}\label{elnuevo}
Any f. p. $k(\n{n})$-module $\um{M}$ decomposes in the following
way
\begin{equation*}
\um{M}\simeq\um{A}_{\lambda_\um{M}}\oplus\left(\bigoplus_{i=1}^n
\uf{F}^i_*\um{B}_{\mu^i_\um{M}}\right)
\oplus\left(\bigoplus_{i=1}^n\uf{F}^i_*\um{C}_{\nu^i_\um{M}}\right)\oplus\uf{M}\uf{S}\um{M}.
\end{equation*}
\end{thm}

This theorem follows from (\ref{metiendo}), (\ref{iny}),
(\ref{tampoco}), (\ref{nuse}), and the following

\begin{lem}\label{buf}
Given two f. p. $k(\n{1})$-modules $\um{M}$ and $\um{N}$ there
exists a $k(\n{n})$-module isomorphism
$\uf{F}^i_*\um{M}\simeq\uf{F}^j_*\um{N}$ $(1\leq i,j\leq n)$ if
and only if one of the following conditions is satisfied:
\begin{itemize}
\item $i=j$ and $\um{M}\simeq\um{N}$,

\item $\um{M}\simeq\um{N}\simeq\um{A}_n$ for some $n\in\NN$.
\end{itemize}
\end{lem}

\begin{proof}
The implication $\Rightarrow$ follows from (\ref{metiendo}),
(\ref{imagenes}) and (\ref{esiso}). On the other hand if
$\um{M}\simeq\um{N}$ then obviously
$\uf{F}^i_*\um{M}\simeq\uf{F}^i_*\um{N}$. Moreover,
$\uf{F}^i_*\um{A}_n$ and $\uf{F}^j_*\um{A}_n$ are isomorphic
$(n\in\NN)$ because both modules are isomorphic to a free
$\ol{T}_n$-controlled $k$-vector space whose basis is a set with
$n$ elements, compare (\ref{ya}).
\end{proof}

As a consequence of Proposition \ref{esigual} and Theorem
\ref{jander} we obtain the next

\begin{cor}\label{unared}
The algebra $k(\n{n})$ has the same representation type as the
$n$-subspace quiver.
\end{cor}

Now Theorem \ref{rep} follows from (\ref{reduc}), (\ref{sesabe})
and (\ref{unared}).



\medskip

We obtain from (\ref{metiendo}), (\ref{imagenes}), (\ref{jander})
and (\ref{elnuevo}) the following presentation of the monoid
$\operatorname{Iso}(\fp( k(\n{n})))$.

\begin{cor}[Classification of f. p. $k(\n{n})$-modules]\label{clasclas}
Let $\set{\ul{V}^{(n,j)}}_{j\in J_n}$ be the set of indecomposable
rigid $n$-subspaces. There is a solution to the decomposition
problem in the category of f. p. $k(\n{n})$-modules given by the
following $1+5n+\card J_n$ elementary modules $(1\leq i\leq n,j\in
J_n)$
$$\uf{F}_*^1\um{A},\;\uf{F}^i_*\um{R},\;\uf{F}^i_*\um{B},\;\uf{F}^i_*\um{B}_\infty,
\;\uf{F}^i_*\um{C},\;\uf{F}^i_*\um{C}_\infty,\;\uf{M}\ul{V}^{(n,j)},\;\;\;$$
and $6n$ elementary isomorphisms $(1\leq i\leq n)$
$$\uf{F}_*^1\um{A}\oplus\uf{F}^i_*\um{R}\simeq\uf{F}^i_*\um{R},\;
\uf{F}^i_*\um{R}\oplus\uf{F}^i_*\um{R}\simeq\uf{F}^i_*\um{R},\;
\uf{F}_*^i\um{B}\oplus\uf{F}^i_*\um{B}_\infty\simeq\uf{F}^i_*\um{B}_\infty,\;$$
$$\uf{F}_*^i\um{B}_\infty\oplus\uf{F}^i_*\um{B}_\infty\simeq\uf{F}^i_*\um{B}_\infty,\;
\uf{F}_*^i\um{C}\oplus\uf{F}^i_*\um{C}_\infty\simeq\uf{F}^i_*\um{C}_\infty,\;
\uf{F}_*^i\um{C}_\infty\oplus\uf{F}^i_*\um{C}_\infty\simeq\uf{F}^i_*\um{C}_\infty.$$
\end{cor}

This classification theorem together with (\ref{losn}) and
(\ref{los4}) complete the proof of Theorem \ref{generel}.

\appendix
\section{Some computations of $\ext_{k(\n{n})}^1$ groups}\label{computo}

The aim of this appendix is to provide with tools and techniques
to compute the $\ext_{k(\n{n})}^1$ group of any pair of f. p.
$k(\n{n})$-modules. This group is in fact a $k$-vector space, so
it is determined by its dimension. Higher $\ext_{k(\n{n})}^*$
groups vanish over f. p. $k(\n{n})$-modules by (\ref{pd}). Since
the functor $\ext_{k(\n{n})}^1$ is biadditive we just have to
compute it over pairs of elementary f. p. $k(\n{n})$-modules, see
(\ref{clasclas}). We shall not make all these computations here
for an arbitrary $n$, but just for $n=1$,
$k(\n{1})=\mathrm{RCFM}(k)$. In addition we show for any $n\in\N$
that the $\ext_{k(\n{n})}^1$ of pairs of f. p. $k(\n{n})$-modules
coming from finite-dimensional $n$-subspaces via the functor
$\uf{M}$ in (\ref{amod}) coincide with their $\ext_{kQ_n}^1$ as
modules over the path algebra. This last vector space is much
easier to compute, since one can use the integral bilinear form of
the quiver $Q_n$, see \cite{taiqf}.



Let $\um{R}$ be the $k$-algebra $\mathrm{RCFM}(k)$ as in Section
\ref{rcfm}. Given two elementary $\mathcal{R}$-modules
$\mathcal{R}/\mat{Y}\mathcal{R}$ $(\mat{Y}\neq 0)$ and
$\mathcal{R}/\mat{Z}\mathcal{R}$, see (\ref{clasifR}), one can
check by using (\ref{semeten}) and basic homological algebra that
there is an isomorphism of $k$-vector spaces
\begin{equation}\label{ext}
\ext^1_\mathcal{R}(\mathcal{R}/\mat{Y}\mathcal{R},\mathcal{R}/\mat{Z}\mathcal{R})\simeq
\frac{\mathcal{R}}{\mathcal{R}\mat{Y}+\mat{Z}\mathcal{R}}.
\end{equation} This formula also holds for $\mat{Z}=0$,
moreover, in this case it is a left-$\mathcal{R}$-module
isomorphism.

\begin{lem}\label{ideal}
We have the following identities
\begin{enumerate}
\item $\mat{A}\um{R}=\set{\mat{R}\in\um{R}\,;\,\mat{r}_{0j}=0
\text{ for all } j\in\NN}$,

\item $(\mat{I}-\mat{A})\um{R}=\set{\mat{R}\in\um{R}\,;\,
\sum_{i\in\NN}\mat{r}_{ij}=0\text{ for all } j\in\NN}$,

\item $(\mat{I}-\mat{A}^t)\um{R}=\set{\mat{R}\in\um{R}\,;\,
\text{given any }i\in\NN,\;\sum_{n\geq i}\mat{r}_{nj}=0\text{ for
almost all }j\in\NN}$,


\item $(\mat{I}-\mat{B}^t)\um{R}=\set{\mat{R}\in\um{R}\,;\,
\text{given }m\in\NN\text{ and }i\leq m,\; \sum_{n\geq
m}\mat{r}_{i+\frac{n(n+1)}{2},j}=0\text{ for almost all
}j\in\NN}$.
\end{enumerate}
\end{lem}

\begin{proof}
One can check that the right-hand-side sets of the statement are
ideals, hence in order to establish the inclusions $\subset$ it is
enough to prove that the matrix defining each left-hand-side ideal
belongs to the corresponding right-hand-side set. This can be
checked by a maybe tedious but straightforward computation.
Suppose now that $\mat{R}$ is a matrix in the right-hand-side set
of (1), (2), (3), or (4) 
, then one can check that the matrix $\mat{C}^1$, $\mat{C}^2$,
$\mat{C}^3$ or $\mat{C}^4$ 
defined as $c^1_{0j}=0$, $\mat{c}^1_{i+1,j}=\mat{r}_{ij}$,
$\mat{c}^2_{ij}=\sum_{n=0}^i\mat{r}_{nj}$,
$\mat{c}^3_{ij}=\sum_{n\geq i}\mat{r}_{nj}$ $(i,j\in\N)$,
$\mat{c}^4_{i+\frac{m(m+1)}{2},j}=\sum_{n\geq m
}\mat{r}_{i+\frac{n(n+1)}{2},j}$ $(i,j\in\NN, m\geq i)$, belongs
to $\um{R}$ and satisfies $\mat{A}\mat{C}^1=\mat{R}$,
$(\mat{I}-\mat{A})\mat{C}^2=\mat{R}$,
$(\mat{I}-\mat{A}^t)\mat{C}^3=\mat{R}$
or $(\mat{I}-\mat{B}^t)\mat{C}^4=\mat{R}$, provided we are in case
(1), (2), (3) or (4). 
Hence we are done.
\end{proof}

By using 
Lemma \ref{ideal} and the involution of the $k$-vector space
$\mathcal{R}$ given by transposition of matrices, one readily
checks that

\begin{lem}\label{ideal2}
We have the following equalities
\begin{enumerate}
\item
$\mathcal{R}(\mat{I}-\mat{A})=\set{\mat{R}\in\mathcal{R}\,;\,
\text{given any }j\in\NN,\;\sum_{n\geq j}\mat{r}_{in}=0\text{ for
almost all }i\in\NN}$,

\item
$\mathcal{R}(\mat{I}-\mat{A}^t)=\set{\mat{R}\in\mathcal{R}\,;\,
\sum_{j\in\NN}\mat{r}_{ij}=0\text{ for all } i\in\NN}$,

\item
$\mathcal{R}(\mat{I}-\mat{B})=\set{\mat{R}\in\mathcal{R}\,;\,
\text{given }m\in\NN\text{ and }j\leq m,\; \sum_{n\geq
m}\mat{r}_{i,j+\frac{n(n+1)}{2}}=0\text{ for almost all
}i\in\NN}$.

\end{enumerate}
\end{lem}

\begin{prop}\label{isos2}
We have that
\begin{enumerate}
\item $\dim
\ext^1_\mathcal{R}(\mathcal{B},\mathcal{R})=2^{\aleph_0}$,

\item $\dim
\ext^1_\mathcal{R}(\mathcal{C},\mathcal{R})=2^{\aleph_0}$,

\item $\dim
\ext^1_\mathcal{R}(\mathcal{B}_\infty,\mathcal{R})=2^{\aleph_0}$,

\item $\dim
\ext^1_\mathcal{R}(\mathcal{C}_\infty,\mathcal{R})=2^{\aleph_0}$,

\item $\dim \ext^1_\mathcal{R}(\mathcal{B},\mathcal{A})=0$,

\item $\dim \ext^1_\mathcal{R}(\mathcal{C},\mathcal{A})=1$,

\item $\dim \ext^1_\mathcal{R}(\mathcal{B}_\infty,\mathcal{A})=0$,

\item $\dim
\ext^1_\mathcal{R}(\mathcal{C}_\infty,\mathcal{A})=\aleph_0$.
\end{enumerate} 
\end{prop}

\begin{proof}

One can check by using (\ref{ideal2}) (1), (\ref{ideal}) (1) and
(\ref{ideal2}) (2), and (\ref{ideal}) (1) and (\ref{ideal2}) (3),
that there are $k$-vector space isomorphisms
$$\frac{\mathcal{R}}{\mathcal{R}(\mat{I}-\mat{A})}\simeq \frac{\prod_{i\in\NN}k}
{\bigoplus_{i\in\NN}k}\;,\;\;\mat{R}+\mathcal{R}(\mat{I}-\mat{A})\mapsto
\left(\sum_{j\in\NN}\mat{r}_{ij}\right)_{i\in\NN}+\bigoplus_{i\in\NN}k,$$
$$\frac{\mathcal{R}}{\mathcal{R}(\mat{I}-\mat{A}^t)+\mat{A}\mathcal{R}}\simeq k\;,\;\;
\mat{R}+(\mathcal{R}(\mat{I}-\mat{A}^t)+\mat{A}\mathcal{R})\mapsto\sum_{j\in\NN}\mat{r}_{0j},$$
$$\frac{\mathcal{R}}{\mathcal{R}(\mat{I}-\mat{B}^t)+\mat{A}\mathcal{R}}\simeq \bigoplus_{j\in\NN} k\;,\;\;
\mat{R}+(\mathcal{R}(\mat{I}-\mat{B}^t)+\mat{A}\mathcal{R})\mapsto\left(\sum_{n\geq
j}\mat{r}_{0,j+\frac{n(n+1)}{2}}\right)_{j\in\NN},$$ hence (1),
(6) and (8) follow from (\ref{ext}).

For any pair of elementary f. p. $\mathcal{R}$-modules the
inequality $\dim \ext^1_\mathcal{R}\leq 2^{\aleph_0}$ follows from
(\ref{ext}) and (\ref{dim}). Now (2) is a consequence of (1) and
(\ref{short}) (1), (3) follows from (1) and the fact that
$\mathcal{B}$ is a direct summand of $\mathcal{B}_\infty$, see
(\ref{clasifR}), and (4) is a consequence of (3) and (\ref{short})
(2).

Given any matrix $\mat{R}\in\mathcal{R}$, if
$\mat{R}^1,\mat{R}^2\in\mathcal{R}$ are the matrices defined by
$\mat{r}^1_{0j}=\mat{r}_{0j}$, $\mat{r}^2_{ij}=\mat{r}_{ij}$
$(i>0)$, and $\mat{r}^n_{ij}=0$ otherwise, then
$\mat{R}=\mat{R}^1+\mat{R}^2$,
$\mat{R}^1\in\mathcal{R}(\mat{I}-\mat{B})$ and
$\mat{R}^2\in\mat{A}\mathcal{R}$, by (\ref{ideal}) (1) and
(\ref{ideal2}) (3), hence
$\mathcal{R}=\mathcal{R}(\mat{I}-\mat{B})+\mat{A}\mathcal{R}$ and
(7) follows by (\ref{ext}). Moreover, since $\mathcal{B}$ is a
direct summand of $\mathcal{B}_\infty$ by (\ref{clasifR}) then (5)
also follows.
\end{proof}

By (\ref{clasifR}), (\ref{medioin}) and (\ref{medioin2}) the
$\ext^1_\mathcal{R}$ group of any other pair of elementary f. p.
$\um{R}$-modules is zero, hence the first extension groups are now
completely computed for f. p. $\um{R}$-modules.

\begin{prop}
For any pair of finite-dimensional $n$-subspaces there is a
natural isomorphism
$$\ext^1_{kQ_n}(\ul{V},\ul{W})\simeq\ext^1_{k(\n{n})}(\uf{M}\ul{V},\uf{M}\ul{W}).$$
\end{prop}

\begin{proof}
Projective representations of $Q_n$ are (arbitrary) direct sums of
the following $n+1$ indecomposable $n$-subspaces,
$$\uf{F}^1(0\rightarrow k),\; \uf{F}^i(k\rightarrow k)\;\;(1\leq i\leq n).$$
Since $\uf{M}$ is an exact full inclusion of categories and any
finite-dimensional representation of $Q_n$ admits a
length-one projective resolution by finite-dimensional projective
representations it is enough to check that
\begin{enumerate}
\item $\ext_{k(\n{n})}^1(\uf{M}\uf{F}^1(0\rightarrow k),
\uf{M}\uf{F}^1(0\rightarrow k))=0$,

\item $\ext_{k(\n{n})}^1(\uf{M}\uf{F}^1(0\rightarrow k),
\uf{M}\uf{F}^i(k\rightarrow k))=0$ $(1\leq i\leq n)$,

\item $\ext_{k(\n{n})}^1(\uf{M}\uf{F}^i(k\rightarrow k),
\uf{M}\uf{F}^1(0\rightarrow k))=0$ $(1\leq i\leq n)$,

\item $\ext_{k(\n{n})}^1(\uf{M}\uf{F}^i(k\rightarrow k),
\uf{M}\uf{F}^j(k\rightarrow k))=0$ $(1\leq i,j\leq n)$.
\end{enumerate}
The resolution constructed in the proof of (\ref{esfp}) shows that
$\uf{M}\uf{F}^1(0\rightarrow k)$ is a projective $k(\n{n})$-module
isomorphic to a 1-dimensional $\ol{T}_n$-controlled $k$-vector
space, hence (1) and (2) hold. Moreover, one can easily check (3)
and (4) by using the definition of $\uf{M}$ in (\ref{amod}) and
the resolutions of $\uf{M}\uf{F}^i(k\rightarrow k)=0$ $(1\leq
i\leq n)$ in the proof of (\ref{esfp}).
\end{proof}

\bibliographystyle{amsplain}
\bibliography{Fernando}
\end{document}